\documentclass[11pt]{amsart}
\usepackage{amsmath}
\usepackage{amsxtra}
\usepackage{amscd}
\usepackage{amsthm}
\usepackage{amsfonts}
\usepackage{amssymb}
\usepackage{eucal}

\numberwithin{equation}{section}

\newtheorem{lem}{Lemma}
\newtheorem{prop}{Proposition}

\theoremstyle{definition}

\theoremstyle{definition}
\newtheorem{thm}{Theorem}

\theoremstyle{remark}
\newtheorem{rem}{Remark}

\numberwithin{thm}{section}
\numberwithin{prop}{section}
\numberwithin{lem}{section}
\numberwithin{rem}{section}

\begin{document}

\newcommand{\thmref}[1]{Theorem~\ref{#1}}
\newcommand{\secref}[1]{Sect.~\ref{#1}}
\newcommand{\lemref}[1]{Lemma~\ref{#1}}
\newcommand{\propref}[1]{Proposition~\ref{#1}}
\newcommand{\corref}[1]{Corollary~\ref{#1}}
\newcommand{\remref}[1]{Remark~\ref{#1}}
\newcommand{\nc}{\newcommand}
\nc{\on}{\operatorname}
\nc{\ch}{\mbox{ch}}
\nc{\Z}{{\mathbb Z}}
\nc{\C}{{\mathbb C}}
\nc{\cond}{|\,}
\nc{\bib}{\bibitem}
\nc{\pone}{\Pro^1}
\nc{\pa}{\partial}
\nc{\F}{{\mathcal F}}
\nc{\arr}{\rightarrow}
\nc{\larr}{\longrightarrow}
\nc{\al}{\alpha}
\nc{\ran}{\rangle}
\nc{\lan}{\langle}
\nc{\W}{{\mathcal W}}
\nc{\gam}{\ovl{\gamma}}
\nc{\Q}{\ovl{Q}}
\nc{\q}{\widetilde{Q}}
\nc{\la}{\lambda}
\nc{\ep}{\epsilon}
\nc{\su}{\widehat{{\mathfrak sl}}_2}
\nc{\gb}{\ovl{{\mathfrak g}}}
\nc{\g}{\overline{{\mathfrak g}}}
\nc{\hh}{\ovl{{\mathfrak h}}}
\nc{\h}{{\mathfrak h}}
\nc{\n}{{\mathfrak n}}
\nc{\ab}{{\mathfrak a}}
\nc{\f}{\widehat{{\mathcal F}}}
\nc{\is}{{\mathbf i}}
\nc{\V}{{\mathcal V}}
\nc{\M}{\widetilde{M}}
\nc{\js}{{\mathbf j}}
\nc{\bi}{\bibitem}
\nc{\laa}{\ovl{\lambda}}
\nc{\fl}{B_-\backslash G}
\nc{\De}{\rho^\vee}
\nc{\G}{{\mathfrak g}}
\nc{\GG}{\widetilde{{\mathfrak g}}}
\nc{\Li}{{\mathcal L}}
\nc{\fp}[2]{\frac{\pa}{\pa u_{#1}^{(#2)}}}
\nc{\Ve}{{\mathcal Vect}}
\nc{\sw}{{\mathfrak s}{\mathfrak l}}
\nc{\La}{\Lambda}
\nc{\ds}{\displaystyle}
\nc{\HH}{\widetilde{\h}}
\nc{\bb}{{\mathfrak b}}
\nc{\wt}{\widetilde}
\nc{\ovl}{\overline}
\nc{\el}{\ell}
\nc{\N}{\ovl{N}_+}
\nc{\nn}{\ovl{\n}}
\nc{\bo}{\ovl{{\mathfrak b}}}
\nc{\pp}{\ovl{p}_{-1}}
\nc{\UU}{{\mathbf U}}
\nc{\VV}{{\mathbf V}}
\nc{\Ss}{{\mathbf S}}
\nc{\QQ}{{\mathbf Q}}
\nc{\vv}{{\mathfrak v}}
\nc{\mc}{\mathcal}
\nc{\ga}{\gamma}
\nc{\mb}{\mathbf}
\nc{\R}{\mathbb R}
\nc{\bu}{\bullet}
\nc{\wh}{\widehat}
\nc{\Gg}{\wh{\G}}
\nc{\AB}{\wh{\ab}}
\nc{\ti}{[t,t^{-1}]}

\title{Five Lectures on Soliton Equations}

\author{Edward Frenkel}

\address{Department of Mathematics, University of California, Berkeley, CA
94720}

\date{July 1997; Revised: November 1997 \\ Contribution to {\em
Surveys in Differential Geometry}, Vol. 3, International Press}

\maketitle

\section*{Introduction}
In these lectures we review a new approach to soliton equations of KdV
type developed by the author together with B.~Feigin and B.~Enriquez
\cite{FF:laws,FF:kdv,EF}.

The KdV equation is the partial differential equation
\begin{equation}    \label{kdv}
v_\tau = 6 v v_z + v_{zzz}.
\end{equation}
Introduced by Korteweg and de Vries in 1895, it has a long and
intriguing history, see, e.g., \cite{Ne}. We will view this equation
as a flow, parametrized by the variable $\tau$, on the space of
functions $v(z)$ in the variable $z$. Perhaps, the most remarkable
aspect of the KdV equation is its {\em complete integrability}, i.e.,
the existence of infinitely many other flows defined by equations of
the type
\begin{equation}    \label{higherkdv}
\pa_{\tau_n} = P_n[v,v_z,v_{zz},\ldots], \quad \quad n=3,4,\ldots,
\end{equation}
which commute with the flow defined by \eqref{kdv} (we set $\tau_1=z$ and
$\tau_2 = \tau$). The important fact is that the right hand sides $P_n$ of
equations \eqref{higherkdv} are polynomials in $v$ and its derivatives,
i.e., {\em differential polynomials in $v$}.

This fact allows us to treat the KdV hierarchy \eqref{higherkdv} in
the following way. Consider the ring $R = \C[v^{(n)}]_{n\geq 0}$ of
polynomials in the variables $v^{(n)}$ (we shall write $v$ for
$v^{(0)}$). Let $\pa_z$ be the derivation of $R$ defined by
the formula $\pa_z \cdot v^{(n)} = v^{(n+1)}$ (so that $v^{(n)} =
\pa_z^n v$). Note that any derivation on $R$ is uniquely determined by
its values on $v^{(n)}$'s via the Leibnitz rule.

Any equation of the form $\pa_\tau v = P[v,v_z,\ldots]$ gives rise to an
evolutionary derivation $D$ of $R$ (i.e., such that commutes with $\pa_z$)
defined by the formula $D \cdot v = P \in R$. The condition of
commutativity with $\pa_z$ implies that $D \cdot v^{(n)} = \pa_z^n \cdot P$, so
that $D$ can be written as
$$D = \sum_{n\geq 0} (\pa_z^n \cdot P) \frac{\pa}{\pa v^{(n)}}.$$ From
this point of view, KdV hierarchy is an infinite set of mutually
commuting evolutionary derivations of $R$. This is the way we will
think of the KdV hierarchy and other, more general, hierarchies
throughout these lectures.

We will not discuss the analytic issues related to the KdV equation,
but will focus on its universal aspects, such as the origin of its
integrability. This approach to soliton equations, initiated by
I.M.~Gelfand and L.A.~Dickey \cite{gd1}, can be summarized as
follows. A particular choice of functional space $F$, in which we want
to look for solutions (we assume that it is an algebra closed under
the action of the derivative $\pa_z$) gives us a specific realization
of the hierarchy. On the other hand, $R$ is nothing but the ring of
functions on the space of $\infty$--jets of functions at a
point. Hence a solution ${\mathbf v}(z,\tau)$ in the realization $F$
gives us a family of homomorphisms $r(t): R \arr F$, which send
$v^{(n)}$ to $\pa_z^n {\mathbf v}(z,t)$. Thus, $R$ can be considered
as a scheme-like object associated with the KdV hierarchy, while
various solutions correspond to various ``points'' of this
``scheme''. This way $R$ captures the universal, realization
independent properties of the KdV equations.

In these notes we will try to answer the question: where do the
infinitely many commuting derivations of KdV hierarchy (and other
hierarchies of KdV type) come from? From the outset, they look rather
mysterious. But in fact there is a very simple explanation. Very
briefly, the answer is as follows:

\begin{enumerate}

\item[(i)] we can identify the vector space with coordinates $v^{(n)},
n\geq 0$, with a homogeneous space of an infinite-dimensional
unipotent algebraic group (a subgroup of a loop group);

\item[(ii)] this homogeneous space has an obvious action of an
infinite-dimensional abelian Lie algebra by vector fields;

\item[(iii)] rewritten in terms of the coordinates $v^{(n)}$, these
vector fields become the KdV flows \eqref{higherkdv}.

\end{enumerate}

Thus we will show that the KdV hierarchy is encoded in the
differential geometry of a loop group.

For technical reasons, this program is easier to realize for a close
relative of the KdV hierarchy -- the so-called modified KdV (or mKdV)
hierarchy; such a hierarchy is associated to an arbitrary affine
Kac-Moody algebra. In Lecture 1 we explain (i)--(iii) for the
generalized mKdV hierarchies. The reader who is only interested in the
main ideas may read just this lecture.

In Lecture 2 we study the mKdV hierarchy in more detail, and write the
equations of the hierarchy explicitly in the so-called zero-curvature
form. In Lecture 3 we compare our approach to mKdV with the earlier
approaches of Drinfeld-Sokolov and Wilson. We remark that in
Sects.~2.1-2.2 and Sects.~3.4-3.6 we follow closely \cite{EF}.

In Lecture 4 we consider the generalized KdV hierarchies from the
point of view of our approach and that of the Drinfeld-Sokolov
reduction \cite{DS}. Some of the results of this lecture (e.g.,
\thmref{vi}) are new. Finally, in Lecture 5 we discuss the Toda field
theories and their connection to the mKdV and KdV hierarchies,
following \cite{FF:kdv}. In particular, we show that the KdV and mKdV
hierarchy are hamiltonian, and their hamiltonians are the integrals of
motion of the corresponding affine Toda field theory.

The original motivation for the formalism explained in these lectures
came from the study of deformations of conformal field theories, where
one needs to show that the integrals of motion of affine Toda field
theories can be quantized (see \cite{Zam,FF:toda}). It turns out that,
in contrast to other approaches, our formalism is well suited for
tackling the quantization problem. In particular, using the results of
Lecture 5 and certain results on quantum groups, we were able to prove
that the Toda integrals of motion can be quantized
\cite{FF:toda,FF:laws}. We briefly discuss this at the end of Lecture
5. Recently these quantum integrals found some spectacular
applications in two-dimensional quantum field theory \cite{BLZ}. We
hope that quantization of our formalism will help elucidate further
the connections between differential geometry of integrable systems
and quantum field theory.

\bigskip

\noindent{\bf Acknowledgements.} I thank S.-T.~Yau for encouraging me
to write this review. I am grateful to B.~Enriquez and B.~Feigin for
their collaboration on the subject of these notes. The notes can serve
as a summary of a lecture course on soliton theory that I gave at
Harvard University in the Spring of 1996. I thank D.~Ben-Zvi for
letting me use his notes of my lectures, and for useful discussions.

I am indebted to B.~Kostant for kindly supplying a proof of
\propref{bk}, which allowed me to simplify the proof of \thmref{vi} on
the structure of the KdV jet space. I also thank L.~Feher for pointing
out an inaccuracy in \cite{EF}, which is corrected here.

The author's research was supported by grants from the Packard
Foundation and the NSF.

\setcounter{section}{1}

\section*{Lecture 1}

\subsection{Background material on affine algebras}    \label{back}

Let $\GG$ be an affine Kac-Moody algebra. It has generators
$e_i,f_i,\al_i^\vee, i=0,\ldots,\el$, and $d$, which satisfy the
standard relations, see \cite{Kac}. The Lie algebra $\GG$ carries a
non-degenerate invariant inner product $(\cdot,\cdot)$. One associates
to $\GG$ the labels $a_i, a_i^\vee, i=0,\ldots,\ell$, the {\em
exponents} $d_i, i=1,\ldots,\ell$, and the {\em Coxeter number} $h$,
see \cite{Ko,Kac}. We denote by $I$ the set of all positive integers,
which are congruent to the exponents of $\GG$ modulo $h$ (with
multiplicities).

The exponents and the Coxeter numbers are given in a table below. Note
that for all affine algebras except $D_{2n}^{(1)}$, each exponent
occurs exactly once. In the case of $D_{2n}^{(1)}$, the exponent
$2n-1$ has multiplicity $2$.

\bigskip

\begin{center}
\renewcommand{\arraystretch}{1.3}
\begin{tabular}{|l|l|l|}
\hline
Type & Coxeter number & \hspace*{20mm}Exponents\\
\hline\hline
$A^{(1)}_n$ & $n+1$ & $1, 2, \ldots, n$ \\
\hline
$A^{(2)}_{2n}$ & $4n+2$ & $1,3,5,\ldots,2n-1,2n+3,\ldots,4n+1$\\
\hline
$A^{(2)}_{2n-1}$ & $4n-2$ & $1,3,5,\ldots,4n-3$\\
\hline
$B^{(1)}_n$ & $2n$ & $1,3,5,\ldots,2n-1$\\
\hline
$C^{(1)}_n$ & $2n$ & $1,3,5,\ldots,2n-1$\\
\hline
$D^{(1)}_n$ & $2n-2$ & $1,3,5,\ldots,2n-3,n-1$\\
\hline
$D^{(2)}_{n+1}$ & $2n+2$ & $1,3,5,\ldots,2n+1$\\
\hline
$D^{(3)}_4$ & $12$ & $1,5,7,11$ \\
\hline
$E^{(1)}_6$ & $12$ & $1,4,5,7,8,11$\\
\hline
$E^{(2)}_6$ & $18$ & $1,5,7,11,13,17$\\
\hline
$E^{(1)}_7$ & $18$ & $1,5,7,9,11,13,17$\\
\hline
$E^{(1)}_8$ & $30$ & $1,7,11,13,17,19,23,29$\\
\hline
$F^{(1)}_4$ & $12$ & $1,5,7,11$\\
\hline
$G^{(1)}_2$ & $6$  & $1,5$\\
\hline
\end{tabular}
\end{center}

\bigskip

The elements $e_i, i=0,\ldots,\el$, and $f_i, i=0,\ldots,\el$,
generate the {\em nilpotent subalgebras} $\n_+$ and $\n_-$ of $\GG$,
respectively. The elements $\al_i^\vee, i=0,\ldots,\el$, and $d$
generate the {\em Cartan subalgebra} $\HH$ of $\GG$. We have:
$\GG=\n_+ \oplus \bb_-$, where $\bb_- = \n_- \oplus \HH$.

The element $$\ds C = \sum_{i=0}^l a_i^\vee \al_i^\vee$$ of $\HH$ is a
{\em central element} of $\GG$. Let $\G$ be the quotient of
$[\GG,\GG]$ by $\C C$. We can identify $\GG$ with the direct sum $\G
\oplus \C C \oplus \C d$. The Lie algebra $\G$ has a Cartan
decomposition $\G = \n_+ \oplus \h \oplus \n_-$, where $\h$ is spanned
by $\al_i^\vee, i=1,\ldots,\el$. Each $x \in \G$ can be uniquely
written as $x_+ + x_-$, where $x_+ \in \n_+$ and $x_- \in \bb_-$.

There exists a unique element $\De \in \HH$, such that $(\al_i,\De)=1,
\forall i=0,\ldots,\el$, and $(d,\De)=0$. But $\HH^*$ is isomorphic to
$\HH$ via the non-degenerate inner product $(\cdot,\cdot)$. We shall
use the same notation for the image of $\De$ in $\HH$ under this
isomorphism. Then $\De$ satisfies: $[\De,e_i]=e_i, [\De,h_i]=0,
[\De,f_i] = -f_i, i=0,\ldots,\el$. The adjoint action of $\De$ on $\G$
defines the {\em principal $\Z$--gradation} on $\G$.

Denote by $\g$ the simple Lie algebra of rank $\el$, generated by
$e_i,f_i,\al_i^\vee, i=1,\ldots,\el$ (the Dynkin diagram of $\g$ is
obtained from that of $\GG$ by removing the $0$th node). It has the
Cartan decomposition $\g = \ovl{\n}_+ \oplus \h \oplus \ovl{\n}_-$,
where $\h$ is the Cartan subalgebra of $\g$, spanned by $\al_i^\vee,
i=1,\ldots,\el$.

\begin{rem}
Affine algebra $\GG$ can be identified with the universal central
extension of the loop algebra $\G^0\ti$, where $\G^0$ is a simple
Lie algebra, or its subalgebra \cite{Kac}. More precisely, if $\GG$ is
non-twisted, then $\G^0=\g$ and $\G = \g\ti$, so that $\GG$ is the
universal central extension of $\g\ti \oplus \C d$:
$$0 \arr \C K \arr \GG \arr \g\ti \oplus \C d\arr 0.$$

A twisted affine algebra $\GG$ can be identified with a subalgebra in
the universal central extension of $\G^0\ti \oplus \C d$.\qed
\end{rem}

\subsection{The principal abelian subalgebra}

Set $$p_{-1} = \sum_{i=0}^\el \frac{(\al_i,\al_i)}{2} f_i,$$ where
$\al_i$'s are the simple roots of $\G$, considered as elements of $\h$
using the inner product. Let $\ab$ be the centralizer of $p_{-1}$ in
$\G$. Recall that we have a principal gradation on $\GG$, such that
$\deg e_i = - \deg f_i = 1, \deg \al_i^\vee = \deg d = 0$. Since
$p_{-1}$ is a homogeneous element with respect to this gradation, the
Lie algebra $\ab$ is a $\Z$--graded subalgebra of $\G$. In particular,
it can be decomposed as $\ab = \ab_+ \oplus \ab_-$, where $\ab_+ = \ab
\cap \n_+$, and $\ab_- = \ab \cap \bb_-$.

\begin{prop}[\cite{Kac1}, Prop.~3.8]    \label{cyclic}
(1) The Lie algebra $\ab$ is abelian.

(2) The homogeneous component of $\ab_\pm$ of degree $i$ with respect
to the principal gradation has dimension equal to the number of
occurencies of $i$ in the set $\pm I$.
\end{prop}

We call $\ab$ the principal abelian subalgebra of $\G$. Its pull-back
to $[\GG,\GG]$ is called the principal Heisenberg subalgebra.

\propref{cyclic} means that for all $\GG$ except $D_{2n}^{(1)}$, each
homogeneous component of $\ab_\pm$ is either zero-dimensional or
one-dimensional, and it is one-dimensional if and only if $i \in
I$. In the latter case we choose a generator $p_i$ of this
component. In the case $\GG = D_{2n}^{(1)}$, the homogeneous component
of degree $i$ congruent to $2n-1$ modulo the Coxeter number $4n-2$ is
two-dimensional. We choose two generators of $\ab$, $p_i^1$ and
$p_i^2$, that span this component.

In particular, we choose
$$\ds p_1 = \sum_{i=0}^l a_i e_i.$$

The properties of $\ab$ that are most important to us are described in
the following fundamental proposition.

\begin{prop}[\cite{Kac1}, Prop. 3.8]    \label{kac}

(1) $\on{Ker} (\on{ad} p_{-1}) = \ab$;

(2) The Lie algebra $\G$ has a decomposition $$\G = \ab \oplus \on{Im}
(\on{ad} p_{-1});$$

(3) With respect to the principal gradation, $\on{Im} (\on{ad} p_{-1})
= \oplus_{j \in \Z} \G_j,$ where $\dim \G_j = \el$, and the map $\on{ad}
p_{-1}: \G_j \arr \G_{j-1}$ is an isomorphism.
\end{prop}

Propositions 1.1 and 1.2 can be derived from B.~Kostant's results
about cyclic elements in simple Lie algebras (see Sect.~6 of
\cite{Ko}).

\begin{rem} L.~Feher has pointed out to us that Prop. 5 of \cite{EF}
was stated incorrectly: it is true only if $n$ is relatively prime to
$h$ (see, e.g., \cite{DF}). However, the case $n=1$ (given above) is
sufficient for the purposes of \cite{EF}, as well as for us here.\qed
\end{rem}

\noindent{\em Examples.} Consider first the case of $$L\sw_2 = \left\{
\begin{pmatrix} a(t) & b(t) \\ c(t) & d(t) \end{pmatrix}, a(t) + d(t) = 0
\right\}.$$ Then
$$\ab = \on{span} \{ p_i \}_{i \on{odd}},$$ where
$$p_{2j+1} = \begin{pmatrix} 0 & t^j \\ t^{j+1} & 0 \end{pmatrix}.$$

More generally, for $\G = L \sw_N$, $$\ab = \on{span} \{ p_i \}_{i
\not{\equiv} 0 \on{mod} N},$$ where $$p_i = (p_1)^i, \quad \quad p_1 =
\begin{pmatrix} 0 & 1 & 0 & \hdots & 0 \\ 0 & 0 & 1 & \hdots & 0 \\
\hdotsfor{5} \\ 0 & 0 & 0 & \hdots & 1 \\ t & 0 & 0 & \hdots & 0
\end{pmatrix}.$$

\begin{rem}    \label{kacpet}
It is known that in the non-twisted case maximal abelian subalgebras
of $\G$ are parametrized by the conjugacy classes of the Weyl group of
$\g$, see \cite{KP}. In particular, $\h\ti$ corresponds to the
identity class, and $\ab$ corresponds to the class of the {\em Coxeter
element}.\qed
\end{rem}

\subsection{The unipotent subgroup}

Let $\wh{\n}_+$ be the completion of the Lie subalgebra $\n_+$ of $\G$
defined as the inverse limit of finite-dimensional nilpotent Lie
algebras $\n_+^{(m)}, m>0$, where $\n_+^{(m)} = \n_+/(\G^0 \otimes t^m
\C[[t]] \cap \n_+)$. In the non-twisted case, $\n_+ = \nn_+ \oplus \g
\otimes t\C[t]$, and so $\wh{\n}_+ = \nn_+ \oplus \g \otimes
t\C((t))$. In general, $\wh{\n}_+$ is a Lie subalgebra of $\G^0((t))$.

Let $N_+^{(m)}$ be the unipotent algebraic group over $\C$
corresponding to $\n_+^{(m)}$. Clearly, the exponential map
$\n_+^{(m)} \arr N_+^{(m)}$ is an isomorphism. We define the
prounipotent proalgebraic group $N_+$ as the inverse limit of
$N_+^{(m)}, m>0$. We have the exponential map $\wh{\n}_+ \arr N_+$,
which is an isomorphism of proalgebraic varieties. Thus, the ring
$\C[N_+]$ of regular functions on $N_+$, is isomorphic
(non-canonically) to the ring of polynomials in infinitely many
variables. Below we will construct such an isomorphism explicitly,
i.e., we will construct explicitly a natural coordinate system on
$N_+$.

\subsection{The action of $\G$}

The group $N_+$ acts on itself from the left and from the right:
$$
n \cdot_L g = n^{-1} g, \quad \quad n \cdot_R g = gn.
$$
In addition, it has a remarkable structure, which is the key in
soliton theory: $N_+$ is equipped with an infinitesimal action of the
Lie algebra $\G$ by vector fields. Let $\ovl{G}$ be the
simply-connected simple algebraic group over $\C$ with Lie algebra
$\g$. In the non-twisted case, there exists an ind-group $G$, whose
set of $\C$--points is $\ovl{G}((t))$. Thus, $G$ can be viewed as the
Lie group of $\G$.  In the twisted case $G$ is defined in a similar
way. The group $G$ has subgroups $N_+$ and $B_-$, where the latter is
the ind-group corresponding to the Lie subalgebra $\bb_-$. We have an
open embedding $N_+ \arr \fl$. Hence the obvious infinitesimal action
of $\G$ on $\fl$ from the right can be restricted to $N_+$.

It is easy to write down the resulting action of $\G$ on $N_+$
explicitly. We can faithfully represent $N_+$ in the Lie algebra of
matrices, whose entries are Taylor power series; for instance, $N_+$
is faithfully represented on $\G^0[[t]]$. Each Fourier coefficient of
such a series defines an algebraic function on $N_+$, and the ring
$\C[N_+]$ is generated by these functions. Hence any derivation of
$\C[N_+]$ is uniquely determined by its action on these functions. We
can write the latter concisely as follows: $\nu \cdot x = y$, where
$x$ is a ``test'' matrix representing an element of $N_+$, i.e. its
$(i,j)$th entry is considered as a regular function $f_{ij}$ on $N_+$,
and $y$ is another matrix, whose $(i,j)$th entry is $\nu \cdot f_{ij}
\in \C[N_+]$.

\begin{rem} One can also take a slightly different point of
view: all formulas below make sense for an arbitrary representation of
$N_+$. Hence instead of picking a particular faithful representation
of $N_+$, one can consider all representations simultaneously. The
corresponding matrix elements generate $\C[N_+]$.\qed
\end{rem}

For $a \in \G$, let us denote by $a^R$ the derivation of $\C[N_+]$
corresponding to the right infinitesimal action of $a$. Also, for $b \in
\n_+$ we denote by $b^L$ the derivation of $\C[N_+]$ corresponding to the
left infinitesimal action of $b$.

In order to simplify notation, we will write below $a x a^{-1}$ instead of
$\on{Ad}_a(x)$.

\begin{lem}
\begin{align}    \label{actionr}
a^R \cdot x &= (x ax^{-1})_+ x, \quad \quad \forall a \in \g,\\
\label{actionl} b^L \cdot x &= -bx, \quad \quad \forall b \in \n_+.
\end{align}
\end{lem}

\begin{proof} Consider a one-parameter subgroup $a(\ep)$ of $G$, such that
$a(\ep) = 1 + \ep a + o(\ep)$. We have: $x \cdot a(\ep) = x + \ep x a
+ o(\ep)$. For infinitesimally small $\ep$ we can factor $x \cdot a(\ep)$
into a product $y_- y_+$, where $y_+ = x + \ep y_+^{(1)} + o(\ep) \in
N_+$ and $y_- = 1 + \ep y_-^{(1)} \in B_-$. We then find that $y_-^{(1)}
x + y_+^{(1)} = x a$, from which we conclude that $y_+^{(1)} = (x
ax^{-1})_+ x$. This proves formula \eqref{actionr}. Formula
\eqref{actionl} is obvious.
\end{proof}

\subsection{The main homogeneous space}

Now we introduce our main object. Denote by $\AB_+$ the completion of
$\ab_+$ in $\wh{\n}_+$. Let $A_+$ be the abelian subgroup of $N_+$,
which is the image of $\AB_+$ in $N_+$ under the exponential
map. Consider the homogeneous space $N_+/A_+$. Since $\G$ acts on
$N_+$ infinitesimally from the right, the normalizer of $\AB_+$ in
$\G$ acts on $N_+/A_+$ infinitesimally from the right. In particular,
$\ab_-$ acts on $N_+/A_+$, and each $p_{-n}, n\in I$, gives rise to a
derivation of $\C[N_+/A_+]$, which we denote by $p_{-n}^R$.

Our goal is to show that the derivations $p_{-n}^R$ are
``responsible'' for the equations of the mKdV hierarchy.

\subsection{The space of jets.}    \label{jets}

Now we change the subject and consider the space $\UU$ of
$\infty$--jets of an $\h$--valued smooth function ${\mb u}(z)$ at
$z=0$. Equivalently, this is the space of functions ${\mathbf u}(z):
{\mathcal D} \arr \h$, where ${\mathcal D}$ is the formal disk,
${\mathcal D} = \on{Spec} \C[[z]]$.

We can view ${\mathbf u}(z)$ as a vector $(u_1,\ldots,u_\el)$, where
$u_i = (\al_i,{\mathbf u}(z))$. The space $\UU$ is therefore the
inverse limit of the finite-dimensional vector spaces $$\UU_m =
\on{span} \{ u_i^{(n)} \}_{i=1,\ldots,\el; n=1,\ldots,m},$$ where $u_i
= (\al_i,{\mathbf u}(0))$, and $u_i^{(n)}=\pa_z^n u_i$ (so, the value
of $u_i^{(n)}$ on ${\mathbf u}(z)$ is given by $\pa_z^n
u_i(0)$). Thus, the ring $\C[\UU]$ of regular functions on $\UU$ is
$\C[u_i^{(n)}]_{i=1,\ldots,\el; n\geq 0}$. The derivative $\pa_z$
gives rise to a derivation $\pa_z$ of $\C[\UU]$, such that $\pa_z
\cdot u_i^{(n)} = u_i^{(n+1)}$. We will write $u_i$ for $u_i^{(0)}$.

Introduce a $\Z$--gradation on $\C[\UU]$ by setting $\deg u_i^{(n)} =
n+1$.

\subsection{Main result}

\begin{thm}[\cite{FF:kdv},Prop.~4]    \label{iso}
{\em There is an isomorphism of rings $$\C[N_+/A_+] \simeq \C[\UU],$$ under
which $p_{-1}^R$ gets identified with $\pa_z$.}
\end{thm}

\begin{proof} We will explicitly construct a homomorphism $\C[\UU] \arr
\C[N_+/A_+]$ and then show that it is actually an isomorphism. In order to
do that, we have to construct regular functions on $N_+/A_+$ corresponding
to $u_i^{(n)} \in \C[\UU], i=1,\ldots,\el; n\geq 0$.

Define a function $\wt{u}_i$ on $N_+$ by the formula
\begin{equation}    \label{wtui}
\wt{u}_i(K) = (\al_i,K p_{-1} K^{-1}), \quad \quad K \in N_+.
\end{equation}
This function is invariant under the right action of $A_+$. Indeed, if
$y \in A_+$, then $y p_{-1} y^{-1} = p_{-1}$, and so $\wt{u}_i(Ky) =
\wt{u}_i(K)$. Hence $\wt{u}_i$ can be viewed as a regular function on
$N_+/A_+$. Next, we define the functions $$\wt{u}_i^{(n)} =
(p_{-1}^R)^n \cdot \wt{u}_i, \quad \quad n>0,$$ on $N_+/A_+$. Consider
the homomorphism $\C[\UU] \arr \C[N_+/A_+]$, which sends $u_i^{(n)}$
to $\wt{u}_i^{(n)}$.

To prove that this homomorphism is injective, we have to show that the
functions $\wt{u}_i^{(n)}$ are algebraically independent. We will do
that by showing that the values of their differentials at the identity
coset $\bar{1} \in N_+/A_+$, $d \wt{u}_i|_{\bar{1}}$, are linearly
independent. Those are elements of the cotangent space to $N_+/A_+$ at
$\bar{1}$, which is canonically isomorphic to
$(\wh{\n}_+/\AB_+)^*$. Using \propref{kac} and the invariant inner
product on $\GG$, we identify $(\wh{\n}_+/\AB_+)^*$ with $(\ab_+)^\perp
\cap \n_-$. By \propref{kac}, with respect to the principal gradation,
$(\ab_+)^\perp = \oplus_{j>0} (\ab_+)^\perp_{-j}$, where $\dim
(\ab_+)^\perp_{-j} = \el$.

Let us first show that the vectors $d \wt{u}_i|_{\bar{1}},
i=1,\ldots,\el$, form a linear basis in $(\ab_+)^\perp_{-1}$. Indeed,
the tangent space to $\bar{1}$ is isomorphic to $\wh{\n}_+/\AB_+$, and
it is clear that the projections of $e_i, i=1,\ldots,\el$, onto
$\wh{\n}_+/\AB_+$ are linearly independent. Hence it suffices to check
that the matrix $[\lan e_i,d \wt{u}_i|_{\bar{1}} \ran]_{1\leq i,j\leq
\el}$ is non-degenerate. But
$$\lan e_i,d \wt{u}_i|_{\bar{1}} \ran = (e_i^L \cdot \wt{u}_j)(\bar{1}),$$
and we find
\begin{align}
(e_i^L \cdot \wt{u}_j)(K) &= - (\al_j,[e_i,K p_{-1} K^{-1}]) \\ &=
- (\al_j,[e_i,p_{-1}] + \ldots) \notag \\ &= - (\al_j,\al_i) \label{neg}
\end{align}
(the dots above stand for the terms lying in $\wh{\n}_+$, and so their
inner product with $\al_j$ is $0$). Thus, we see that the matrix
$[\lan e_i,d \wt{u}_i|_{\bar{1}} \ran]_{1\leq i,j\leq \el}$ coincides
with minus the symmetrized Cartan matrix of $\g$, $-
[(\al_i,\al_j)]_{1\leq i,j\leq \el}$, and hence is
non-degenerate. Therefore the vectors $d \wt{u}_i|_{\bar{1}},
i=1,\ldots,\el$, are linearly independent.

Now, by definition, $$d \wt{u}^{(m+1)}_i|_{\bar{1}} = \on{ad} p_{-1}
\cdot d \wt{u}^{(m)}_i|_{\bar{1}}.$$ But according to \propref{kac},
$\on{ad} p_{-1}: (\ab_+)^\perp_{-j} \arr (\ab_+)^\perp_{-j-1}$ is an
isomorphism for all $j>0$. Hence the vectors $d
\wt{u}^{(m)}_i|_{\bar{1}}$ are all linearly independent. Therefore the
functions $\wt{u}_i^{(n)}, i=1,\ldots,\el, n\geq 0$, are algebraically
independent and our homomorphism $\C[\UU] \arr \break \C[N_+/A_+]$ is
an embedding.

Now let us compare gradations. Consider the derivation $(\rho^\vee)^R$
on $N_+$. It is clear that for any $v \in \ab_+, [\rho^\vee,v] \in
\ab_+$ (see \secref{back}). Hence $(\rho^\vee)^R$ is a well-defined
derivation on $\C[N_+/A_+]$. We take $-(\rho^\vee)^R$ as the gradation
operator on $\C[N_+/A_+]$ (to simplify notation, we will denote it by
$-\rho^\vee$).

We have by formula \eqref{actionr},
\begin{align*}
(-\rho^\vee \cdot \wt{u}_i)(K) &= -(\al_i,[(K \De K^{-1})_+,K p_{-1}
K^{-1})] \\ &= -(\al_i,[K \De K^{-1},K p_{-1} K^{-1})] + (\al_i,[(K
\De K^{-1})_-,K p_{-1} K^{-1})] \\ &= - (\al_i,K [\De,p_{-1}] K^{-1})
+ ([\al_i,(K \De K^{-1})_-],K p_{-1} K^{-1})] \\ &= (\al_i,K p_{-1}
K^{-1}) = \wt{u}_i(K).
\end{align*}
Hence the degree of $\wt{u}_i$ equals $1$. Since $[-\De,p_{-1}] =
p_{-1}$, the degree of $p_{-1}^R$ also equals $1$, and so $\deg
wt{u}_i^{(n)} = n+1$.

On the other nand, the function $u_i^{(n)}$ has degree $n+1$ with
respect to the $\Z$--gradation on $\C[\UU]$. Hence our homomorphism
$\C[\UU] \arr \C[N_+/A_+]$ is homogeneous of degree $0$. To show its
surjectivity, it suffices to prove that the characters of $\C[\UU]$
and $\C[N_+/A_+]$ coincide. Here by character of a $\Z$--graded vector
space $V$ we understand the formal power series $$\on{ch} V = \sum_{n
\in \Z} \dim V_n q^n,$$ where $V_n$ is the homogeneous subspace of $V$
of degree $n$. Clearly,
$$\on{ch} \C[\UU] = \prod_{n>0} (1-q^n)^{-\el}.$$

Now, it is well-known that as a $\G$--module, $\C[N_+]$ is isomorphic
to the module $M^*_0$ contragradient to the Verma module with highest
weight $0$ (these modules are defined in Lecture 4). Hence, as an
$\n_+$ module, it is isomorphic to the restricted dual $U(\n_+)^*$ of
$U(\n_+)$. This implies that $\C[N_+/A_+]$ is isomorphic to the space
of $\ab_+$--invariants on $U(\n_+)^*$.

By the Poincare-Birkhoff-Witt theorem, $U(\n_+)^*$ has a basis dual to
the basis of lexicographically ordered monomials in $U(\n_+)$. We
can choose this basis in such a way that each element is the product
$A B$, where $A$ is the product of $p_n, n \in I$ (a basis in
$\ab_+$), and $B$ is a product of elements of some basis in $\on{Im}
(\on{ad} p_{-1}) \cap \n_+$. Then the elements of $U(\n_+)^*$ dual to
the monomials of the second type form a basis of the space of
$\ab_+$--invariants on $U(\n_+)^*$. Hence the character of
$\C[N_+/A_+]$ equals that of a free polynomial algebra with $\el$
generators of each positive degree. Thus its character coincides with
the character of $\C[\UU]$, and this completes the proof.
\end{proof}

\subsection{Definition of the mKdV hierarchy}

Having identified $\C[\UU]$ with $\C[N_+/A_+]$ in the theorem above,
we can now view the derivation $p_{-n}^R, n \in I$ of the latter, as a
derivation of $\C[\UU]$, which we denote by $\pa_n$. Thus, we obtain
an {\em infinite set of commuting derivations} of the ring of
differential polynomials $\C[u_i^{(n)}]$, of degrees equal to the
exponents of $\G$ modulo the Coxeter number.

We call this set the {\em mKdV hierarchy} associated to the affine
algebra $\G$.

According to \thmref{iso}, to each map ${\mathcal D} \arr \h$ we can
attach a map ${\mathcal D} \arr N_+/A_+$ (here ${\mathcal D}$ is the
formal disk). From analytic point of view, this means that we have a
mapping that to a every smooth function ${\mb u}(z)$ from, say,
${\mathbb R}$ to $\h$, assigns a smooth function $K(z): \R \arr
N_+/A_+$. Moreover, this mapping is {\em local} in the sense that
$K(z)$ depends on ${\mb u}$ only through the values of its derivatives
at $z$ (the jets of ${\mb u}(z)$ at $z$). Due to this property, ${\mb
u}(z)$ can just as well be a smooth function on a circle, or an
analytic function on a domain in $\C$.

Next we say that on $N_+/A_+$ we have a family of commuting vector
fields $p_{-n}^R$ that come from the right action of the Lie algebra
$\ab_-$. The flows of the mKdV hierarchy are the corresponding flows
on the space of smooth function ${\mb u}(z): {\mathbb R} \arr
\h$. Note that the function ${\mb u}(z)$ does not really have to be
smooth everywhere; it may have singularities at certain points.

In principle, we have now answered the question posed in the
Introduction as to where commuting derivations acting on rings of
differential polynomials come from. But we would like to have explicit
formulas for these derivations. From the construction itself we know
that the first of them, $\pa_1$, is just $\pa_z$. But we don't know
the rest. On the other hand, the mKdV hierarchies have been previously
defined by other methods, and we want to make sure that our definition
coincides with the earlier definitions.

These issues will be the subject of the next two lectures. We will first
write the derivations $\pa_n$ in the so-called zero-curvature form and then
compare our definition of the mKdV hierarchy with the other definition to
see that they are equivalent.

\setcounter{section}{2}
\setcounter{subsection}{0}
\setcounter{equation}{0}
\setcounter{thm}{0}
\setcounter{prop}{0}
\setcounter{lem}{0}
\setcounter{rem}{0}

\section*{Lecture 2}

\subsection{Zero curvature representation}

Let $\omega_i^\vee \in \h$ be the fundamental coweights that satisfy:
$(\al_i,\omega_j^\vee) = \delta_{i,j}$. Define ${\mathbf u}$ to be the
element
$$\sum_{i=1}^\el \omega_i^\vee \otimes u_i \in \h \otimes \C[\UU].$$

Denote by $\Gg$ the completion of $\G$, which is the inverse limit of
$\G/(\G^0 \otimes t^m\C[t] \cap \G)$. Thus, if $\G=\g\ti$, then
$\Gg=\g((t))$, and in general $\Gg$ is a Lie subalgebra of
$\G^0((t))$.

For each $K \in N_+/A_+$, $K p_{-n}K^{-1}$ is a well-defined element
of $\Gg$. Thus we obtain a algebraic map $N_+/A_+ \arr \Gg$, or,
equivalently, an element of $\Gg \otimes \C[N_+/A_+]$.\footnote{More
precisely, it lies in the completed tensor product, e.g., if $\Gg =
\g((t))$, then by $\Gg \otimes \C[N_+/A_+]$ we mean $(\g \otimes
\C[N_+/A_+])((t))$.} By abuse of notation, we denote this element by
$K p_{-n}K^{-1}$, and its projection on ${\mathfrak b}_- \otimes
\C[N_+/A_+]$ by $(K p_{-n}K^{-1})_-$. Since $\C[N_+/A_+] \simeq
\C[u_i^{(n)}]$, by \thmref{iso}, we can apply to it $1 \otimes \pa_m$,
which we denote simply by $\pa_m$.

More explicitly, recall that the Lie algebra $\Gg$ can be realized as a
Lie subalgebra of the Lie algebra $\G^0((t))$, where $\G^0$ is a
finite-dimensional simple Lie algebra; e.g., in the non-twisted case,
$\G^0 = \g$. If we choose a basis in $\G^0$, we can consider an
element of $\Gg$ as a matrix, whose entries are Laurent power series.
The entries of the matrix $Kp_{-n}K^{-1}$ are Laurent series in $t$
whose coefficients are regular functions on $N_+/A_+$. Hence, under
the isomorphism $\C[N_+/A_+] \simeq \C[u_i^{(n)}]$, each coefficient
corresponds to a differential polynomial in $u_i$'s. Applying $\pa_m$
to $K p_{-n}K^{-1}$ means applying $\pa_n$ to each of these
coefficients.

We are going to
prove the following result.

\begin{thm}[\cite{EF}, Theorem 2]    \label{zeroc}
\begin{equation}    \label{dsmkdv1}
[\pa_z + p_{-1} + {\mathbf u},\pa_n + (Kp_{-n}K^{-1})_-] = 0.
\end{equation}
\end{thm}

The equation \eqref{dsmkdv1} can be rewritten as
\begin{equation}    \label{dsmkdv}
\pa_n {\mathbf u} = \pa_z (Kp_{-n}K^{-1})_- + [p_{-1}+{\mathbf
u},(Kp_{-n}K^{-1})_-].
\end{equation}
As explained above, this equation expresses $\pa_n u_i$ in terms of
differential polynomials in $u_i$'s. Since, by construction, $\pa_n$
commutes with $\pa_1 \equiv \pa_z$, formula \eqref{dsmkdv} uniquely
determines $\pa_n$ as an evolutionary derivation of $\C[u_i^{(n)}]$.
Thus, we obtain an explicit formula for $\pa_n$.

Note that formula \eqref{dsmkdv1} looks like the zero curvature
condition on a connection defined on a two-dimensional space -- hence
the name ``zero curvature representation''.

\begin{rem}
One should be careful in distinguishing the two variables: $t$ and
$z$, and the two factors, $\Gg$ and $\C[N_+/A_+]$, in the formulas
below. The variable $z$ is a ``dynamical variable'' connected with
$\C[N_+/A_+]$ due to its isomorphism with $\C[u_i^{(n)}]$; in
particular, $\pa_z$ is the derivation of $p_{-1}^R$ of
$\C[N_+/A_+]$. On the other hand, $t$ is just the formal variable
entering the definition of the affine algebra $\Gg$. We remark that in
earlier works on soliton theory, $t$ was denoted by $\la^{-1}$, and
$\la$ was called the {\em spectral parameter}.\qed
\end{rem}

\subsection{Proof of \thmref{zeroc}}

We will actually prove a stronger result:

\begin{prop}    \label{zcgen}
For $K \in N_+/A_+$,
\begin{equation}    \label{zcgenf}
[\pa_m + (Kp_{-m}K^{-1})_-,\pa_n + (Kp_{-n}K^{-1})_-] = 0, \quad \quad
\forall m,n \in I.
\end{equation}
\end{prop}

Let us first explain how to derive \thmref{zeroc} from
\propref{zcgen}. We need to specialize \eqref{zcgenf} to $m=1$ and to
determine $(Kp_{-1}K^{-1})_-$ explicitly.

\begin{lem}    \label{p-1}
\begin{equation}    \label{boldu}
(Kp_{-1}K^{-1})_- = p_{-1} + {\mathbf u}.
\end{equation}
\end{lem}

\begin{proof}
It is clear that $(Kp_{-1}K^{-1})_- = p_{-1} + x$, where $x \in \h$. Hence
we need to show that $x={\mathbf u}$, or, equivalently, that $(\al_i,x) =
u_i, i=1,\ldots,\el$. We can rewrite the latter formula as
$u_i=(\al_i,(Kp_{-1}K^{-1})_-)$, and hence as
$u_i=(\al_i,Kp_{-1}K^{-1})$. But this is exactly the definition of the
function $\wt{u}_i$ on $N_+/A_+$, which is the image of $u_i \in \C[\UU]$
under the isomorphism $\C[\UU] \arr \C[N_+/A_+]$. 
\end{proof}

Now specializing $m=1$ in formula \eqref{zcgenf} and using \lemref{p-1} we
obtain formula \eqref{dsmkdv1}.

In order to prove formula \eqref{zcgenf}, we need to find an explicit
formula for the action of $\pa_n$ on $Kp_{-m}K^{-1}$.

It follows from formula \eqref{actionr} that
\begin{equation}    \label{actionc}
a^R \cdot x  v x^{-1} = [(x ax^{-1})_+,x  v x^{-1}], \quad
\quad a,v \in \G.
\end{equation}
If $a$ and $v$ are both elements of $\ab$, then formula \eqref{actionc}
does not change if we multiply $x$ from the right by an element of
$A_+$. Denote by $K$ the coset of $x$ in $N_+/A_+$. Then we can write:
\begin{equation}    \label{actionc1}
\pa_n \cdot K v K^{-1} = [(Kp_{-n}K^{-1})_+,K v K^{-1}], \quad \quad v
\in \ab.
\end{equation}

\vspace*{5mm}
\noindent {\em Proof of \propref{zcgen}.} Substituting $v=p_{-m}$ into
formula \eqref{actionc1}, we obtain:
$$\pa_n \cdot K p_{-m} K^{-1} = [(Kp_{-n}K^{-1})_+,K p_{-m} K^{-1}].$$
Hence $$\pa_n \cdot (K p_{-m} K^{-1})_- = [(Kp_{-n}K^{-1})_+,K p_{-m}
K^{-1}]_- = [(Kp_{-n}K^{-1})_+,(K p_{-m} K^{-1})_-]_-.$$ Therefore we
obtain:
\begin{align*}
& [\pa_m + (Kp_{-m}K^{-1})_-,\pa_n + (Kp_{-n}K^{-1})_-] \\ = & \pa_m \cdot
(K p_{-n}K^{-1})_- - \pa_n \cdot (K p_{-n} K^{-1})_- + [(K p_{-m}
K^{-1})_-,(K p_{-n} K^{-1})_-] \\ = & [(Kp_{-m}K^{-1})_+,(K p_{-n}
K^{-1})_-]_- - [(Kp_{-n}K^{-1})_+,K p_{-m} K^{-1}]_- \\ + & [(K p_{-m}
K^{-1})_-,(K p_{-n} K^{-1})_-]_-.
\end{align*}

Adding up the first and the last terms, we obtain $$[K p_{-m}
K^{-1},(Kp_{-n}K^{-1})_-]_- - [(Kp_{-n}K^{-1})_+,K p_{-m} K^{-1}]_- =$$
$$= [K p_{-m} K^{-1},K p_{-n} K^{-1}]_- = 0,$$
and \propref{zcgen} is proved.\qed

\subsection{Recurrence relation}    \label{rr}

Now we have at our disposal explicit formulas for equations of the mKdV
hierarchy. But from the practical point of view, it is still difficult to
compute explicitly terms like $(K p_{-n} K^{-1})_-$, since the coset $K$
has been defined in Lecture 1 in a rather abstract way. In this section we
will exhibit a simple property of $K$, which will enable us to compute
everything via a straightforward recursive algorithm.

Given a vector space $W$, we will write $W[\UU]$ for $W \otimes
\C[\UU]$.

Recall that for any $v \in \ab$, $K v K^{-1}$ denotes the element of
$\Gg[\UU] = \Gg \otimes \C[\UU]$, determined by the embedding $N_+/A_+
\arr \Gg$, which maps $K$ to $K v K^{-1}$.

\begin{lem}    \label{centr}
\begin{equation}    \label{comm}
[\pa_z + p_{-1} + {\mathbf u},K v K^{-1}]=0, \quad \quad \forall v \in
\ab.
\end{equation}
\end{lem}

\begin{proof} Using formula \eqref{actionc1} and \lemref{p-1} we obtain:
\begin{align*}
\pa_z (K v K^{-1}) &= [(Kp_{-1}K^{-1})_+,K v K^{-1}] \\ &= -
[(Kp_{-1}K^{-1})_-,K v K^{-1}] \\ &= - [p_{-1} + {\mathbf u},K v K^{-1}].
\end{align*}
\end{proof}

Now suppose that $\V \in \Gg[\UU]$ satisfies
\begin{equation}    \label{recu}
[\pa_z + p_{-1} + {\mathbf u}(z),\V]=0.
\end{equation}
Then we can decompose $\V$ with respect to the principal gradation on
$\Gg$ (the principal gradation on $\Gg$ should not be confused with the
gradation on $\C[\UU] = \C[N_+/A_+]$!):
$$\V = \sum_{m \geq -n} \V_m, \quad \quad \V_{-n} \neq 0,$$ and write
equation \eqref{recu} in components. Recall that $\deg p_{-1} = -1, \deg
{\mathbf u} = 0$, and $\deg \pa_z = 0$.

The first equation that we obtain, in degree $-n-1$, reads:
$$[p_{-1},\V_{-n}] = 0.$$ By \propref{kac}, it implies that $\V_{-n}
\in \ab$. Hence $\V_{-n} = p_{-n}$ times an element of $\C[\UU]$,
where $n \in I \cup -I$.

Other equations have the form
\begin{equation}    \label{recur}
\pa_z \V_m + [p_{-1},\V_{m-1}] + [{\mathbf u},\V_m] = 0,
\quad \quad m \geq -n.
\end{equation}
Recall from \propref{kac} that $\G = \ab \oplus \on{Im} (\on{ad}
p_{-1})$. So we can split each $x \in \Gg$ into a sum $x = x^0 + x^\perp$,
where $x^0 \in \ab$ and $x^\perp \in \on{Im} (\on{ad} p_{-1})$. Furhermore,
we can split our equations \eqref{recur} into two parts, lying in $\ab$ and
$\on{Im} (\on{ad} p_{-1})$. Then we obtain two equations:
\begin{align*}
\pa_z \V_m^0 + [p_{-1},\V_{m+1}]^0 + [{\mathbf u},\V_m]^0 &= 0, \\
\pa_z \V_m^\perp + [p_{-1},\V_{m+1}]^\perp + [{\mathbf u},\V_m]^\perp &= 0.
\end{align*}
But clearly, $[p_{-1},\V_{m+1}]^0 = 0$, $[{\mathbf u},\V_m^0]^0 = 0$, and
$[p_{-1},\V_{m+1}^0] = 0$. Also, $[p_{-1},\V_{m+1}^\perp]^\perp =
[p_{-1},\V_{m+1}^\perp]$. Hence we obtain:
\begin{align}
\pa_z \V_m^0 + [{\mathbf u},\V_m^\perp]^0 &= 0, \label{odin} \\
\pa_z \V_m^\perp + [p_{-1},\V_{m+1}^\perp] + [{\mathbf u},\V_m]^\perp &= 0
\label{dva}.
\end{align}

We already know that $\V_{-n} \in \ab[\UU]$, and so $\V_{-n}^\perp =
0$. Now we can find $\V_{-n}^0$ and $\V_{-n+1}^\perp$. First we obtain
from equation \eqref{odin}:
$$\pa_z \V_{-n}^0 = - [{\mathbf u},\V_{-n}^\perp] = 0,$$ which implies
that $\V_{-n} = p_{-n}$ times a {\em constant} factor. Without loss of
generality we assume that $\V_{-n} = p_{-n}$. Next, we obtain from
\eqref{dva}, $$[p_{-1},\V_{-n+1}^\perp] = - [{\mathbf
u},\V_{-n}]^\perp,$$ which uniquely determines $\V_{-n+1}^\perp$,
since $\on{ad} p_{-1}$ is invertible on $\on{Im} (\on{ad} p_{-1})$
(see \propref{kac}).

Now assume by induction that we know $\V_k^0, -n \leq k \leq m-1$ and
$\V_k^\perp, -n \leq k \leq m$. Then we obtain from formulas \eqref{odin}
and \eqref{dva}:
\begin{align}
\pa_z \V_m^0 &= - [{\mathbf u},\V_m^\perp]^0, \label{odin1} \\
\V_{m+1}^\perp &= (\on{ad} p_{-1})^{-1} \left( - \pa_z \V_m^\perp -
[{\mathbf u},\V_m]^\perp \right) \label{dva1}.
\end{align}
Consider $[{\mathbf u},\V_m^\perp]^0 \in \ab$. It can only be non-zero
if $m \in \pm I$. In that case it is equal to $p_m \otimes P_m, P_m
\in \C[\UU]$, or $p_m^1 \otimes P_m^1 + p_m^2 \otimes P_m^2$, if the
multiplicity of the exponent is $2$. If each $P_m^i$ is a total
derivative, $P_m^i = \pa_z Q_m^i$, then we find from equation
\eqref{odin1}:
$$\V_m^0 = \sum_i p_m^i \otimes Q_m^i.$$ After that we can solve
equation \eqref{dva1} for $\V_{m+1}^\perp$. But we know that $K p_{-n}
K^{-1} = p_{-n} + \ldots$ does satisfy equation \eqref{recur}. Hence
we know for sure that each $P_m^i$ is a total derivative, and equation
\eqref{odin1} can be resolved (i.e., the antiderivative can be taken
whenever we need it). Therefore we can find all $\V_m$'s following
this algorithm.

Note that the only ambiguity is introduced when we take the
antiderivative, since we then have the freedom of adding $c_m \cdot
p_m$ to $\V_m^0$, where $c_m \in \C$. This simply reflects the fact
that adding $K p_m K^{-1} = p_{-m} + \ldots$ with $m>-n$ to $K p_{-n}
K^{-1}$, we do not violate \eqref{recur}. Thus, we have proved the
following

\begin{prop}    \label{classify}
Let $\V$ be an element of $\Gg[\UU]$ satisfying
\begin{equation}    \label{nov}
[\pa_z + p_{-1} + {\mathbf u},\V] = 0.
\end{equation}
Then $\V = K v K^{-1}$, where $v \in \ab$. Moreover, if $\V = p_{-n}
+$ terms of degree higher than $-n$ with respect to the principal
gradation on $\G$ then $v = p_{-n} +$ terms of degree higher than
$-n$.
\end{prop}

Now recall that each $\V_n$ is an element of $\Gg[\UU]$, i.e., its
matrix entries (in a particular representation of $\Gg$) are Laurent
power series in $t$, with coefficients from the ring $\C[\UU]$. That
ring has a $\Z$--gradation (which should not be confused with the
principal gradation on $\G$ itself!). The corresponding gradation on
$\C[N_+/A_+]$ is defined by the vector field $-\De$. Using this
gradation, for each $n$ we can distinguish a canonical solution $\V_n$
of the recurrence relations \eqref{recur} that equals precisely $K
p_{-n} K^{-1}$.

\begin{lem} Let $x \in \G$ be homogeneous of degree $k$. Then the
function $f \in \C[N_+/A_+]$ given by the formula $f(K) =
(x,Kp_{-n}K^{-1})$ is a homogeneous element of $\C[N_+/A_+]$ of degree
$n-m$.
\end{lem}

\begin{proof} We have to show that $\De \cdot f = (m-n)f$. But
\begin{align*}
(\De \cdot f)(K) &= (x,[(K \De K^{-1})_+,K p_{-n} K^{-1}]) \\ &= (x,[K \De
K^{-1},K p_{-n} K^{-1}]) - (x,[(K \De K^{-1})_-,K p_{-n} K^{-1}]) \\
&= (x,K[\De,p_{-n}]K^{-1}) + ([\De,x],K p_{-n} K^{-1}]) \\
&= (m-n) (x,K p_{-n} K^{-1}) = (m-n) f(K),
\end{align*}
and the lemma is proved.
\end{proof}

Our recurrence relations \eqref{odin} and \eqref{dva} imply that if
$x$ has degree $m$, then $(x,\V_n)$ has a term of degree $n-m$ in
$\C[\UU] \simeq \C[N_+/A_+]$ and possibly other terms of smaller
degrees which result from addition of constants when taking
anti-derivatives at the previous steps of the recursive procedure. It
is clear that there is a unique solution $\V_n$ of equation
\eqref{recu}, such that $(x,\V_n)$ is homogeneous of degree $n-m$, and
this solution is $K p_{-n} K^{-1}$.

Now we have an algorithm for finding $K p_{-n} K^{-1}$. Hence we can
construct explicitly the $n$th equations of the mKdV hierarchy by inserting
$(K p_{-n} K^{-1})_-$ into equation \eqref{dsmkdv1}. In the next section we
will do that in the case of $L\sw_2$.

\subsection{Example of $L\sw_2$}

We have $${\mathbf u} = \begin{pmatrix} \dfrac{1}{2} u & 0 \\
0 & - \dfrac{1}{2} u \end{pmatrix}.$$ Hence
\begin{equation}    \label{Lop}
\pa_z + p_{-1} + {\mathbf u} = \pa_z + \begin{pmatrix} \dfrac{1}{2} u & t^{-1}
\\ 1 & - \dfrac{1}{2} u \end{pmatrix}.
\end{equation}

The first equation of the hierarchy is 
$$[\pa_z + p_{-1} + {\mathbf u},\pa_n + (Kp_{-n}K^{-1})_-] = 0.$$

By formula \eqref{boldu}, $(K p_{-1} K^{-1})_- = p_{-1} + {\mathbf u}$.
Hence the first equation reads
$$[\pa_z + p_{-1} + {\mathbf u},\pa_1 + p_{-1} + {\mathbf u}] = 0,$$
which is equivalent to $$\pa_1 {\mathbf u} = \pa_z {\mathbf u},$$ as
we already know (note that this is true for an arbitrary $\G$).

Now we consider the next equation in the case of $L\sw_2$, corresponding to
$n=3$. To write it down, we have to compute $(K p_{-3} K^{-1})_-$. In order
to do that we compute the first few terms of $K p_{-3} K^{-1}$ following
the algorithm of the previous section.

We introduce a basis of $L\sw_2$ homogeneous with respect to the principal
gradation:
$$p_{2j+1} = \begin{pmatrix} 0 & t^j \\ t^{j+1} & 0 \end{pmatrix}, \quad
\quad q_{2j+1} = \begin{pmatrix} 0 & t^j \\ -t^{j+1} & 0 \end{pmatrix},$$
$$r_{2j} = \begin{pmatrix} t^j & 0 \\ 0 & - t^j \end{pmatrix}.$$ Then we
write $$K p_{-3} K^{-1} = p_{-3} + R_{-2} r_{-2} + P_{-1} p_{-1} + Q_{-1}
q_{-1} + R_0 r_0 + \ldots.$$ The other terms are irrelevant for us now as
we are only interested in $(K p_{-3} K^{-1})_-$.

Equation \eqref{recu} reads
$$[\pa_z + p_{-1} + \frac{1}{2} u r_0,\pa_3 + p_{-3} + R_{-2} r_{-2} + P_{-1}
p_{-1} + Q_{-1} q_{-1} + R_0 r_0 + \ldots] = 0.$$ Rewriting it in
components, we obtain the following equations:
$$2 R_{-2} - u = 0, \quad \quad \pa_z R_{-2} - 2 Q_{-1} = 0,$$
$$\pa_z P_{-1} + u Q_{-1} = 0, \quad \quad \pa_z Q_{-1} + u P_{-1} - 2 R_0.$$
We find from these equations: $$R_{-2} = \frac{1}{2} u, \quad Q_{-1} =
\frac{1}{4} \pa_z u, \quad P_{-1} = - \frac{1}{8} u^2, \quad R_0 = -
\frac{1}{16} u^3 + \frac{1}{8} \pa_z^2 u.$$

The equation of the mKdV hierarchy corresponding to $n=3$ now reads:
$$\pa_3 \begin{pmatrix} \dfrac{1}{2} u & t^{-1} \\ 1 & - \dfrac{1}{2} u
\end{pmatrix} =$$ $$\left[ \begin{pmatrix} \dfrac{1}{2} u t^{-2} - \left(
\dfrac{1}{16} u^3 - \dfrac{1}{8} \pa_z^2 u \right) & t^{-2} + \left(
- \dfrac{1}{8} u^2 + \dfrac{1}{4} \pa_z u \right) t^{-1} \\ t^{-1} - \left(
\dfrac{1}{8} u^2 + \dfrac{1}{4} \pa_z u \right) & - \dfrac{1}{2} u
t^{-2} + \left( \dfrac{1}{16} u^3 - \dfrac{1}{8} \pa_z^2 u \right)
\end{pmatrix} , \pa_z + \begin{pmatrix} \dfrac{1}{2} u & t^{-1} \\ 1 & -
\dfrac{1}{2} u \end{pmatrix} \right].$$ This is equivalent to the equation
\begin{equation}    \label{mkdv1}
\pa_3 u = \frac{3}{8} u^2 \pa_z u - \frac{1}{4} \pa_z^3 u,
\end{equation}
which is the {\em mKdV equation}. It is related to the KdV equation
\eqref{kdv} by a change of variables (see Lecture 4).

\setcounter{section}{3}
\setcounter{subsection}{0}
\setcounter{equation}{0}
\setcounter{thm}{0}
\setcounter{prop}{0}
\setcounter{lem}{0}
\setcounter{rem}{0}

\section*{Lecture 3}

In the previous lecture we have written down explicitly the equations
of the mKdV hierarchy in the zero curvature form. In this lecture we
will discuss another approach to these equations, which is due to
Drinfeld and Sokolov and goes back to Zakharov and Shabat. We will
then establish the equivalence between the two approaches.

\subsection{Generalities on zero curvature equations}

Let us look at the general zero curvature equation
\begin{equation}    \label{zc}
[\pa_z + p_{-1} + {\mathbf u},\pa_{\tau} + V] = 0.
\end{equation}

The element $p_{-1}$ has degree $-1$ with respect to the principal
gradation of $\G$, while ${\mathbf u}$ has degree $0$ (with respect to
the gradation on $\G$). This makes finding an element $V$ that
satisfies \eqref{zc} a non-trivial problem. Indeed, equation
\eqref{zc} can be written as
\begin{equation}    \label{next}
\pa_{\tau} {\mathbf u} = [\pa_z + p_{-1} + {\mathbf u},V].
\end{equation}
The left hand side of \eqref{next} has degree $0$. Therefore $V$ should
be such that the expression in the right hand side of \eqref{next} is
concentrated in degree $0$.

Such elements can be constructed by the following trick (see
\cite{ZS,DS,W}). Suppose we found some ${\mathcal V} \in \Gg[\UU]$
which satisfies
\begin{equation}    \label{next1}
[\pa_z + p_{-1} + {\mathbf u},{\mathcal V}] = 0.
\end{equation}
We can split ${\mathcal V}$ into the sum ${\mathcal V}_{+} + {\mathcal
V}_{-}$ of its components of positive and non-positive degrees with
respect to the principal gradation. Then $V = {\mathcal V}_{-}$ has the
property that the right hand side of \eqref{next} has degree $0$. Indeed,
from \eqref{next1} we find
$$[\pa_z + p_{-1} + {\mathbf u}(z),{\mathcal V}_{-}] = - [\pa_z + p_{-1} + {\mathbf
u}(z),{\mathcal V}_{+}],$$ which means that both commutators have neither
positive nor negative homogeneous components. Therefore equation
\eqref{next} makes sense. Thus, the problem of finding elements $V$ such
that \eqref{next} makes sense, and hence constructing the mKdV hierarchy,
has been reduced to the problem of finding solutions of equation
\eqref{next1}.

In \propref{classify} we described the solutions of equation \eqref{next1}
in terms of the isomorphism $\UU \simeq N_+/A_+$, which attaches $K \in
N_+/A_+$ to ${\mathbf u}$.

There are other ways of solving this problem. In the next section we will
discuss the approach of Drinfeld and Sokolov \cite{DS}, closely related to
the dressing method of Zakharov and Shabat \cite{ZS}. Another approach,
proposed by Wilson \cite{W}, will be mentioned briefly later.

\subsection{Drinfeld-Sokolov approach}

Consider the group $N_+[\UU]$ of $\UU$--points of $N_+$. Thus, in a
representation of $N_+$ elements of $N_+[\UU]$ are matrices, whose
entries are differential polynomials in $u_i$'s.

\begin{prop}[\cite{DS}, Prop.~6.2]    \label{conj}
There exists an element $M \in N_+[\UU]$, such that
\begin{equation}    \label{ds}
M^{-1} \left( \pa_z + p_{-1} + {\mathbf u}(z) \right) M = \pa_z + p_{-1}
+ \sum_{i\in I} h_i p_i,
\end{equation}
where $h_i \in \C[\UU], \forall i \in I$. $M$ is defined uniquely up
to right multiplication by an element of $A_+[\UU]$.
\end{prop}

\begin{proof}
Since the exponential map $\wh{\n}_+ \arr N_+$ is an isomorphism, we can
write $M = \exp {\mathbf m}$, where ${\mathbf m} \in \wh{\n}_+[\UU]$. For
any $\al \in \G$, we have:
$$M^{-1} \al M = \sum_{n \geq 0} \frac{(-1)^n}{n!} (\on{ad} {\mathbf m})^n
\cdot \al.$$ Decompose ${\mathbf m}$ with respect to the principal gradation:
$${\mathbf m} = \sum_{j>0} {\mathbf m}_j.$$ Then for the $n$th component of
equation \eqref{ds} reads
\begin{equation}    \label{n=0}
[p_{-1},{\mathbf m}_1] + {\mathbf u} = 0,
\end{equation}
for $n=0$, and
\begin{equation}    \label{urav}
[p_{-1},{\mathbf m}_{n+1}] + \text{terms involving} \, \, {\mathbf
m}_i, i\leq n, \, \, \text{and} \, \, {\mathbf u} = h_n p_n
\end{equation}
for $n>0$ (here we set $h_n=0$, if $n \not{\!\!\in} I$). Since
${\mathbf u} \in \on{Im} (\on{ad} p_{-1})$, we can solve equation \eqref{n=0}
and start the inductive process. At the $n$th step of induction we would
have found ${\mathbf m}_i, i\leq n$ and $h_i, i<n$. We can now split the
$n$th equation \eqref{urav} into two parts, lying in $\ab$ and $\on{Im}
(\on{ad} p_{-1})$. We can then set ${\mathbf m}_{n+1} = (\on{ad}
p_{-1})^{-1}$ of the $\on{Im} (\on{ad} p_{-1})$ part, and $h_n p_n =$ the
$\ab$ part of the equation (everything evidently works even if the
multiplicity of $n$ is $2$). Clearly, both ${\mathbf m}_n$ and $h_n$ is
local. This inductive procedure gives us a unique solution $M$ of equation
\eqref{ds}, such that each ${\mathbf m}_{n+1}$ lies in $\on{Im} (\on{ad}
p_{-1})$.
\end{proof}

Now observe that formula \eqref{ds} implies:
\begin{equation}    \label{aha}
[\pa_z + p_{-1} + {\mathbf u},M v M^{-1}] = 0, \quad \quad \forall v
\in \ab.
\end{equation}
Hence the zero curvature equations
\begin{equation}    \label{oldds}
[\pa_z + p_{-1} + {\mathbf u},\pa_{\tau_n} + (M p_{-n} M^{-1})_-] = 0
\end{equation}
make sense. Drinfeld and Sokolov called these equations the mKdV hierarchy
associated to $\G$. We will now show that these equations coincide with our
equations \eqref{dsmkdv1}.

\subsection{Equivalence of two constructions}

In both \cite{DS} and \cite{FF:kdv}, one assigns to a jet $(u_i^{(n)})
\in \on{Spec} \C[\UU]$, a coset in $N_+/A_+$. We now show that they
coincide.

\begin{thm} [\cite{EF}, Theorem 3]   \label{ident}
{\em The cosets $M$ and $K$ in $N_+/A_+$ assigned in \cite{DS} and
\cite{FF:kdv}, respectively, to a jet $(u_i^{(n)})$, coincide.}
\end{thm}

\begin{proof}
According to formula \eqref{aha}, $$[\pa_z + p_{-1} + {\mathbf u},M
p_{-1} M^{-1}] =0.$$ Since $M p_{-1} M^{-1} = p_{-1} +$ terms of
degree higher than $-1$ with respect to the principal gradation on
$\G$, we obtain from \lemref{classify} that $M p_{-1} M^{-1} = K v
K^{-1}$, where $v \in \ab$. Let us show that this implies that
$v=p_{-1}$ and that $K=M$ in $N_+/A_+$.

Indeed, from the equality $$Mp_{-1}M^{-1} = K v K^{-1}$$ we obtain
that $(K^{-1} M) p_{-1} (K^{-1} M)^{-1}$ lies in $\ab$. We can
represent $K^{-1} M$ as $\exp y$ for some $y \in \wh{\n}_+$. Then
$(K^{-1} M) p_{-1} (K^{-1} M)^{-1} = v$ can be expressed as a linear
combination of multiple commutators of $y$ and $p_{-1}$: $$e^y p_{-1}
(e^y)^{-1} = \sum_{n\geq 0} \frac{1}{n!}  (\on{ad} y)^n \cdot
p_{-1}.$$ We can write $y = \sum_{j>0} y_j$, where $y_j$ is the
homogeneous component of $y$ of principal degree $j$. It follows from
\propref{kac} that $\wh{\n}_+ = \AB_+ \oplus \on{Im} (\on{ad}
p_{-1})$. Therefore each $y_j$ can be further split into a sum of
$y_j^0 \in \AB_+$ and $y_j^1 \in \on{Im} (\on{ad} p_{-1})$.

Suppose that $y$ does not lie in $\AB_+$. Let $j_0$ be the smallest number
such that $y_{j_0}^1 \neq 0$. Then the term of smallest degree in $e^y
p_{-1} (e^y)^{-1}$ is $[y_{j_0}^1,p_{-1}]$ which lies in $\on{Im} (\on{ad}
p_{-1})$ and is non-zero, because $\on{Ker} (\on{ad} p_{-1}) =
\AB_+$. Hence $e^y p_{-1} (e^y)^{-1}$ can not be an element of $\AB_+$.

Therefore $y \in \AB_+$ and so $K^{-1} M \in A_+$, which means that $K$ and
$M$ represent the same coset in $N_+/A_+$, and that $v=p_{-1}$.
\end{proof}

\begin{rem} In the proof of Theorem 3 of \cite{EF} one has to replace
$p_{-n}$ by $p_{-1}$.\qed
\end{rem}

{}From the analytic point of view, \propref{conj} produces the same
thing as \thmref{iso}: it associates to a smooth function ${\mb u}(z):
\R \arr \h$, a smooth function $M(z): \R \arr N_+/A_+$
(cf. Sect.~1.8). What we have just shown is that the two constructions
are equivalent, and $K(z)=M(z)$.

Having established that $K=M$, we see that the two definitions of mKdV
hierarchy: \eqref{dsmkdv1} and \eqref{oldds}, coincide. Thus, now we
have two equivalent ways of constructing these equations. One of them
is based on \thmref{iso} identifying $\on{Spec} \C[u_i^{(n)}]$ with
$N_+/A_+$, and the other uses the ``dressing operator'' $M$ defined by
formula \eqref{ds}. One can then use an element of $N_+/A_+$ ($K$ or
$M$) to solve equation \eqref{next1}, and produce the mKdV equation
\eqref{next}.

Wilson \cite{W} proposed another approach, in which one solves
\eqref{next1} directly without constructing first the dressing
operator. Recall from \secref{rr} that when one tries to solve equation
\eqref{next1}, one has to be able to take the antiderivative at certain
steps. If we know in advance the existence of solution (as we do in the two
other approaches), then this is insured automatically. Wilson \cite{W} gave
an argument, which demonstrates directly that the antiderivatives can be
taken at each step. Thus, he proved the existence of solution ${\mathcal V}$ of
\eqref{next1} of the form ${\mathcal V} = p_{-n} +$ terms of higher degree. He
then constructed the equations of the mKdV hierarchy using these ${\mathcal
V}$. We now see that all three approaches are equivalent.

In the next two sections we will give another, more illuminating,
explanation of the fact that $K=M$.

\begin{rem}
There exist generalizations of the mKdV hierarchies which are
associated to abelian subalgebras of $\G$ other then $\ab$ (see
\remref{kacpet}). It is known that the Drinfeld-Sokolov approach can
be applied to these generalized hierarchies \cite{DHM,HM}. On the
other hand, our approach can also be applied; in the case of the
non-linear Schr\"{o}dinger hierarchy, which corresponds to the
homogeneous abelian subalgebra of $\G$, this has been done by Feigin
and the author \cite{FF:nls}. The results of this lecture can be
extended to establish the equivalence between the two approaches in
this general context.\qed
\end{rem}

\subsection{Realization of $\C[N_+]$ as a polynomial ring.}

The approach to the mKdV and affine Toda equations used in
\cite{FF:kdv} and here is based on \thmref{iso} which identifies
$\C[N_+/A_+]$ with the ring of differential polynomials
$\C[u_i^{(n)}]_{i=1,\ldots,\el;n\geq 0}$. In this section we add to
the latter ring new variables corresponding to $A_+$ and show that the
larger ring thus obtained is isomorphic to $\C[N_+]$. In this and the
following two sections we follow closely \cite{EF}.

Consider $u_i^{(n)}, i=1,\ldots,\el; n\geq 0$, as $A_+$--invariant
regular functions on $N_+$. Recall that $$u_i(x)=(\al_i,x
p_{-1}x^{-1}), \quad \quad x \in N_+,$$ Now choose an element
$\chi$ of $\HH$, such that $(\chi,C) \neq 0$. Introduce the regular
functions $\chi_n, n\in I$, on $N_+$ by the formula:
\begin{equation}    \label{chin}
\chi_n(x)=(\chi,x p_{-n}x^{-1}), \quad \quad x \in N_+.
\end{equation}

\begin{thm}[\cite{EF}, Theorem 5]    \label{identif}
$\C[N_+] \simeq \C[u_i^{(n)}]_{i=1,\ldots,\el;n\geq 0} \otimes
\C[\chi_n]_{n\in I}.$
\end{thm}

\begin{proof} Let us show that the functions $u_i^{(n)}$'s and $\chi_n$'s
are algebraically independent. In order to do that, let us compute the
values of the differentials of these functions at the origin. Those are
elements of the cotangent space to the origin, which is isomorphic to the
dual space $\wh{\n}_+^*$ of $\wh{\n}_+$.

Using the invariant inner product on $\GG$ we identify $\wh{\n}_+^*$
with $\n_-$. By \propref{kac} $\n_- = \ab_- \oplus
(\ab_+)^\perp$. Recall that with respect to the principal gradation,
$(\ab_+)^\perp = \oplus_{j>0} (\ab_+)^\perp_{-j}$, where $\dim
(\ab_+)^\perp_{-j} = \el$, and $\on{ad} p_{-1}: (\ab_+)^\perp_{-j}
\arr (\ab_+)^\perp_{-j-1}$ is an isomorphism for all $j>0$.

By construction of $u_i$'s given in the proof of \thmref{iso}, $du_i|_1,
i=1,\ldots,\el$, form a basis of $(\ab_+)^\perp_{-1}$, and hence
$du^{(n)}_i|_1, i=1,\ldots,\el$, form a basis of
$(\ab_+)^\perp_{-n-1}$. Thus, the covectors $du^{(n)}_i|_1, i=1,\ldots,\el;
n\geq 0$, are linearly independent. Let us show now that the covectors
$d\chi_n|_1$ are linearly independent from them and among themselves. For
that it is sufficient to show that the pairing between $dF_m|_1$ and $p_n$
is non-zero if and only if $n=m$. But we have:
\begin{equation}    \label{cox}
p_n^R \cdot (\chi,x p_{-m}x^{-1}) = (\chi,x[p_n,p_{-m}]x^{-1}).
\end{equation}
Since $p_n, n \in \pm I$, generate a Heisenberg subalgebra in $\GG$
(see \cite{Kac}, Lemma 14.4),
\begin{equation}    \label{nonzero}
[p_n,p_{-m}] = \al_n C \delta_{n,m},
\end{equation}
where $\al_n \neq 0, \forall n \in I$. Therefore the pairing between
$dF_m|_1$ and $p_n$ is equal to $\al_n (\chi,C) \delta_{n,-m}$, and it
is non-zero if and only if $n=m$.

Thus, the functions $u_i^{(n)}$'s and $\chi_n$'s are algebraically
independent. Hence we have an embedding
$\C[u_i^{(n)}]_{i=1,\ldots,\el;n\geq 0} \otimes \C[\chi_n]_{n\in I}
\arr \C[N_+]$. But the characters of the two spaces with respect to
the principal gradation are both equal to
$$\prod_{n\geq 0} (1-q^n)^{-\el} \prod_{i\in I} (1-q^i)^{-1}.$$ Hence this
embedding is an isomorphism.
\end{proof}

\subsection{Another proof of \thmref{ident}.}

Let $c_n = (\al_n (\chi,C))^{-1}$, where $\al_n, n\in I$, denote the
non-zero numbers determined by formula \eqref{nonzero}.

\begin{prop}    \label{localm}
Let $K$ be an element of $N_+$. We associate to it another element of
$N_+$, $$\ovl{K} = K \exp \left( - \sum_{n\in I} c_n p_{n} \chi_n(K)
\right)$$ In any finite-dimensional representation of $N_+$, $\ovl{K}$ is
represented by a matrix whose entries are Taylor series with coefficients
in the ring of differential polynomials in $u_i, i=1,\ldots,\el$

The map $N_+ \arr N_+$ which sends $K$ to $\ovl{K}$ is constant on the
right $A_+$--cosets, and hence defines a section $N_+/A_+ \arr N_+$.
\end{prop}

\begin{proof}
Each entry of $$\ovl{K} = K \exp \left( - \sum_{n \in I}
c_n p_{n} \chi_n(K) \right)$$ is a function on
$N_+$. According to \thmref{iso} and \thmref{identif}, to prove the
proposition it is sufficient to show that each entry of $\ovl{K}$ is
invariant under the right action of $\ab_+$. By formula \eqref{cox} we
obtain for each $m \in I$:
$$
p_m^R \cdot \chi_n = c_n^{-1} \delta_{n,-m}, 
$$
and hence
$$
p_m^R \cdot \left( K \exp \left( - \sum_{n\in I} c_n p_{n}\chi_n(K)
\right) \right) =$$ $$K p_m \exp \left( - \sum_{n\in I} c_n
p_{n}\chi_n \right) + K \exp \left( - \sum_{n\in I} c_n p_{n}\chi_n
\right) (-p_m) = 0.
$$
Therefore $\ovl{K}$ is right $\ab_+$--invariant.

To prove the second statement, let $a$ be an element of $A_+$ and let
us show that $\overline{Ka} = \ovl{K}$. We can write: $a = \exp
\left(\sum_{i \in I} \al_i p_i \right)$. Then according to formulas
\eqref{chin} and \eqref{cox}, $\chi_n(Ka) = (\chi,Ka p_{-n} a^{-1}
K^{-1}) = (\chi,K p_{-n} K^{-1}) + c_n^{-1} \al_n = \chi_n(K) +
c_n^{-1} \al_n$. Therefore
$$\overline{Ka} = Ka \exp \left( - \sum_{n \in I} \al_n p_n - \sum_{n
\in I} c_n p_n \chi_n(K) \right) = \ovl{K}.$$
\end{proof}

Consider now the matrix $\ovl{K}$. According to \propref{localm}, the
entries of $\ovl{K}$ are Taylor series with coefficients in
differential polynomials in $u_i$'s. In other words, $\ovl{K} \in
N_+[\UU]$. Hence we can apply to $\ovl{K}$ any derivation of
$\C[u_i^{(n)}]$, in particular, $\pa_n = p_{-n}^R$.

\begin{lem}    \label{dsd}
In any finite-dimensional representation of $N_+$, the matrix of $\ovl{K}$
satisfies:
\begin{equation}    \label{ds1}
\ovl{K}^{-1} (\pa_n + (\ovl{K} p_{-n} \ovl{K}^{-1})_-) \ovl{K} = \pa_n +
p_{-n} - \sum_{i \in I} c_i (p_{-n}^R \cdot \chi_i) p_i.
\end{equation}
\end{lem}

\begin{proof}
Using formula \eqref{actionr}, we obtain:
\begin{align*}
& \ovl{K}^{-1} (\pa_n + (\ovl{K} p_{-n} \ovl{K}^{-1})_-) \ovl{K} = \pa_n +
\ovl{K}^{-1} (p_{-n}^R \ovl{K}) + \ovl{K}^{-1} (\ovl{K} p_{-n}
\ovl{K}^{-1})_- \ovl{K} \\ & = \pa_n + \ovl{K}^{-1} (\ovl{K} p_{-n}
\ovl{K}^{-1})_+ \ovl{K} - \sum_{i \in I} c_i (p_{-n}^R \cdot \chi_i) p_i +
\ovl{K}^{-1} (\ovl{K} p_{-n} \ovl{K}^{-1})_- \ovl{K} \\ & = \pa_n + p_{-n}
- \sum_{i \in I} c_i (p_{-n}^R \cdot \chi_i) p_i,
\end{align*}
which coincides with \eqref{ds1}.
\end{proof}

Let now $K$ be the point of $N_+/A_+$ assigned to the jet
$(u_i^{(n)})$ by \thmref{iso}. Then $\ovl{K}$ is a well-defined
element of $N_+$ corresponding to $K$ under the map $N_+/A_+ \arr N_+$
defined in \propref{localm}. By construction, $\ovl{K}$ lies in the
$A_+$--coset $K$.

According to \lemref{p-1}, $(\ovl{K} p_{-1} \ovl{K}^{-1})_- = p_{-1} +
{\mathbf u}$. Letting $n=1$ in formula \eqref{ds1}, we obtain:
$$\ovl{K}^{-1} (\pa_z + p_{-1} + {\mathbf u}(z)) \ovl{K} = \pa_z +
p_{-1} - \sum_{i \in I} c_i (p_{-1}^R \cdot \chi_i) p_i.$$ This shows
that $\ovl{K}$ gives a solution to equation \eqref{ds}, and hence lies
in the $A_+$--coset of the Drinfeld-Sokolov dressing operator
$M$. Therefore the cosets of $K$ and $M$ coincide, and this
completes our second proof of \thmref{identif}.

\subsection{Baker-Akhiezer function}

In this section we adopt the analytic point of view. Recall that in
Lecture 1 we constructed a map, which assigns to every smooth function
${\mb u}(z): {\mathbb R} \arr \h$, a smooth function $K(z): {\mathbb
R} \arr N_+/A_+$, with the property that the action of $n$th mKdV flow
on ${\mb u}(z)$ translates into the action of $p_{-n}^R$ on
$N_+/A_+$. One can ask whether it is possible to lift this map to the
one that assigns to ${\mb u}$ a function ${\mathbb R} \arr N_+$ with
the same property. Our results from the previous section allow us to
lift the function $K(z)$ to the function $\ovl{K}(z): {\mathbb R} \arr
N_+$, using the section $N_+/A_+ \arr N_+$. But one can check easily
that the action of the mKdV flows does not correspond to the action of
$\ab_-^R$ on $\ovl{K}(z)$. However, its modification
$$\wt{K}(z) = \ovl{K}(z) \exp \left( \sum_{n\in I} c_n p_{n}
\chi_n(\ovl{K}) \right)$$ does satisfy the desired property.

Unfortunately, $\chi_n \not{\!\!\in} \C[u_i^{(n)}] \subset
\C[N_+]$, so $\wt{K}(z)$ can not be written in terms of differential
polynomials in $u_i$ (i.e, the jets of ${\mb u}(z)$), so $\wt{K}(z)$
is ``non-local''. Nevertheless, one can show that $H_n = p_{-1}^R
\cdot \chi_n$ does belong to $\C[u_i^{(n)}]$ (in fact, this $H_n$ can
be taken as the density of the hamiltonian of the $n$th mKdV equation,
see Lecture 5). Therefore we can write:
$$\wt{K}(z) = \ovl{K}(z) \exp \left( \sum_{n\in I} c_n p_{n}
\int_{-\infty}^z H_n dz \right).$$

Now let ${\mb u}({\mb t})$, where ${\mb t} = \{ t_i \}_{i \in
I}$ and $t_i$'s are the times of the mKdV hierarchy, be a solution of
the mKdV hierarchy. The Baker-Akhiezer function $\Psi({\mb t})$
associated to ${\mb u}({\mb t})$, is by definition a solution of
the system of equations
\begin{equation}    \label{bakergen}
(\pa_n + (Kp_{-n}K^{-1})_-) \Psi = 0, \quad \quad \forall n \in I.
\end{equation}
In particular, for $n=1$ we have:
\begin{equation}    \label{baker}
(\pa_z + p_{-1} + {\mb u}(z)) \Psi = 0.
\end{equation}

{}From formula \eqref{ds1} we obtain the following explicit formula
for the solution of the system \eqref{bakergen} with the initial
condition is $\Psi({\mb 0}) = \wt{K}({\mb 0})$:
\begin{equation*}
\Psi({\mb t}) = \wt{K}({\mb t}) \exp \left( - \sum_{i \in I}
p_{-i}t_i \right).
\end{equation*}
But by construction, $$\wt{K}({\mb t}) = \left( \wt{K}({\mb 0})
\Gamma({\mb t}) \right)_+,$$ where $$\Gamma({\mb t}) = \exp \left(
\sum_{i \in I} p_{-i}t_i \right)$$ and $g_+$ denotes the projection of
$g \in B_- \cdot N_+ \subset G$ on $N_+$ (it is well-defined for
almost all $t_i$'s). Hence we obtain:
$$\Psi({\mb t}) = \left( \wt{K}({\mb 0}) \Gamma({\mb t}) \right)_+
\Gamma({\mb t})^{-1}.$$

Similar formula for the Baker-Akhiezer function has been obtained by
G.~Segal and G.~Wilson \cite{SW,Wil}. Following the works of the Kyoto
School \cite{DJKM}, they showed that Baker-Akhiezer functions
associated to solutions of the KdV equations naturally ``live'' in the
Sato Grassmannian, and that the flows of the hierarchy become linear
in these terms (see also \cite{Ch}).

We have come to the same conclusion in a different way. We have
identified the mKdV variables directly with coordinates on the big
cell of $B_- \backslash G/A_+$, and constructed a map ${\mb u}(z) \arr
K(z)$. It is then straightforward to check, as we did above, that the
Baker-Akhiezer function is simply a lift of this map to $B_-
\backslash G$. In the case of KdV hierarchy, we obtain a map to
$\ovl{G}[t^{-1}] \backslash G$ (see the next lecture), which for
$\ovl{G}=SL_n$ is the formal version of the Grassmannian that
Segal-Wilson had considered.

\setcounter{section}{4}
\setcounter{subsection}{0}
\setcounter{equation}{0}
\setcounter{thm}{0}
\setcounter{prop}{0}
\setcounter{lem}{0}
\setcounter{rem}{0}

\section*{Lecture 4}

In this lecture we define the generalized KdV hierarchy associated to
an affine algebra $\G$. We then show the equivalence between our definition
and that of Drinfeld-Sokolov.

Throughout this lecture we restrict ourselves to non-twisted affine algebras.

\subsection{The left action of $\N$}

Let $\ovl{\n}_+$ be the Lie subalgebra of $\n_+$ generated by $e_i,
i=1,\ldots,\el$. Thus, $\ovl{\n}_+$ is the upper nilpotent subalgebra of
the simple Lie algebra $\g$ -- the ``constant subalgebra'' of $\GG$, whose
Dynkin diagram is obtained by removing the $0$th nod of the Dynkin digram
of $\GG$. Let $\N$ be the corresponding subgroup of $N_+$. The group $\N$
acts on the main homogeneous space $N_+/A_+$ from the left.

\begin{lem}
The action of $\N$ on $N_+/A_+$ is free.
\end{lem}

\begin{proof} We need to show that for each $K \in N_+/A_+$, the map $\tau$ from
$\ovl{\n}_+$ to the tangent space of $N_+/A_+$ at each point $K$ is
injective. But the tangent space at $K$ is naturally isomorphic to $\wh{\n}_+/(K
\ab K^{-1})$, and $\tau$ is the compositon of the embedding $\ovl{\n}_+
\arr \wh{\n}_+$ and the projection $\wh{\n}_+ \arr \wh{\n}_+/(K \ab K^{-1})$. Thus, we
need to show that $\ovl{\n}_+ \cap K \ab K^{-1} = 0$ for any $K \in
N_+/A_+$. This is obvious because each element of $\ovl{\n}_+$ is a
constant element of $\G$, i.e. does not depend on $t$, while each element
of $K \ab K^{-1}$ does have a $t$--dependence.
\end{proof}

\subsection{The KdV jet space}

Now we show that $\C[\N \backslash N_+/A_+]$ can be identified with a ring
of differential polynomials.

We define the functions $v_i: N_+ \arr \C, i=1,\ldots,\el$, by the
formula:
\begin{equation}    \label{vik}
v_i(K) = (f_0,K p_{-d_i} K^{-1}), \quad \quad K \in
N_+.
\end{equation}
It is clear that these function are right $A_+$--invariant. We also
find:
$$
(e_k^L \cdot v_i)(K) = (f_0,[e_k,K p_{-n} K^{-1}]) = ([f_0,e_k],K
p_{-n} K^{-1}) = 0
$$
for all $k=1,\ldots,\el$. Therefore $v_i$ is left $\nn_+$--invariant,
and hence left $\N$--invariant. Thus, each $v_i$ gives rise to a
regular function on $\N \backslash N_+/A_+$. Let us compute the degree
of $v_i$. We obtain:
\begin{align*}
(-\rho^\vee \cdot v_i)(K) &= -(f_0,[(K \De K^{-1})_+,K p_{-d_i}
K^{-1})] \\ &= -(f_0,[K \De K^{-1},K p_{-d_i}
K^{-1})] + (f_0,[(K \De K^{-1})_+,K p_{-d_i}
K^{-1})] \\ &= -(f_0,K[\De,p_{-1}]K^{-1}) + ([f_0,\De],K p_{-d_i}
K^{-1})] \\ &= (d_i+1) (f_0,K p_{-d_i} K^{-1}) = (d_i+1) v_i(K).
\end{align*}
Hence $v_i$ is homogeneous of degree $d_i+1$.

Now denote $v_i^{(n)} = (p_{-1}^R)^n \cdot v_i$.

\begin{thm}    \label{vi}
$$\C[\N \backslash N_+/A_+] \simeq \C[v_i^{(n)}]_{i=1,\ldots,\el;n\geq
0}.$$
\end{thm}

\begin{proof}
We follow the proof of \thmref{iso}. First we prove the algebraic
independence of the functions $v_i^{(n)}$. To do that, we need to
establish the linear independence of the values of their differentials
$dv_i^{(n)}|_{\bar{1}}$ at the double coset $\bar{1}$ of the identity
element. Note that $dv_i^{(n)}|_{\bar{1}} = (\on{ad} p_{-1})^n \cdot
dv_i|_{\bar{1}}$.

Using the invariant inner product on $\GG$, we identify the cotangent
space to $\bar{1}$ with ${\mc N} = (\ab_+ \oplus \nn_+)^\perp \cap
\n_-$. It is a graded subspace of $\n_-$ with respect to the principal
gradation. Moreover, the action of $\on{ad} p_{-1}$ on ${\mc N}$ is
free. One can choose homogeneous elements $\beta_1,\ldots,\beta_\el
\in {\mc N}$, such that $\{ (\on{ad} p_{-1})^n \beta_i
\}_{i=1,\ldots,\el;n\geq 0}$ is a basis in $(\ab_+ \oplus \nn_+)^\perp
\cap \n_+^*$. It is also easy to see that $\deg
\beta_i=-(d_i+1)$. This follows from \propref{kac} and the fact
\cite{Ko} that the set of degrees of a elements of a homogeneous basis
in $\ovl{\n}_-$ is $\cup_{i=1}^\el \{ 1,\ldots,d_i \}$.

We claim that $dv_i|_{\bar{1}}$ can be taken as the elements
$\beta_i$. Let us compute $dv_i|_{\bar{1}}$. From formula \eqref{vik}
we find, in the same way as in the proof of \thmref{iso}:
$$
(x,dv_i|_{\bar{1}}) = - (f_0,[x,p_{-d_i}]) = ([f_0,p_{-d_i}],x) \quad
\quad \forall x \in \n_+.
$$
Hence
\begin{equation}    \label{dvi}
dv_i|_{\bar{1}} = [f_0,p_{-d_i}]
\end{equation}
as an element of ${\mc N} \subset \n_-$. This formula shows that the
degree of $dv_i|_{\bar{1}}$ equals $d_i+1$ as we already know.

According to \propref{kac}, $(\ab_+)^\perp \subset \n_-$ has an
$\el$--dimensional component $(\ab_+)^\perp_{-j}$ of each negative
degree $-j$ (with respect to the principal gradation). Furthermore,
since $\ovl{\n}_-$ has no homogeneous components of degrees less than
or equal to $-h$, minus the Coxeter number, we obtain that
$(\ab_+)^\perp_{-j} = {\mc N}_{-j}$ for all $j\geq h$. Now let
\begin{equation}    \label{gai}
\ga_i = (\on{ad} p_{-1})^{h-d_i-1} \cdot dv_i|_{\bar{1}} = (\on{ad}
p_{-1})^{h-d_i-1} \cdot [f_0,p_{-d_i}] \in (\ab_+)^\perp_{-h} =
{\mc N}_{-h}
\end{equation}
for all $i=1\ldots,\el$. Recall that the operator $\on{ad} p_{-1}:
(\ab_+)^\perp_{-j} \arr (\ab_+)^\perp_{-j-1}$ is an isomorphism. If
the vectors $\ga_i$ are linearly independent, then the covectors
$$
d v_i^{(j-d_i-1)}|_{\bar{1}} = (\on{ad} p_{-1})^{j-d_i-1} \cdot
[f_0,p_{-d_i}] \in {\mc N}_{-j},
$$
are linearly independent for each $j>1$ (here we ignore $d
v_i^{(n)}|_{\bar{1}}$, if $n<0$). We see that the linear independence
of the covectors $\{ dv_i^{(n)}|_{\bar{1}} \}_{i=1,\ldots,\el;n\geq
0}$ is equivalent to the following

\begin{prop}    \label{bk}
The vectors $\{ \ga_i \}_{i=1,\ldots,\el}$ are linearly independent.
\end{prop}

The proof is given in the Appendix. This completes the proof of the
algebraic independence of the functions $v_i^{(n)}$. Hence we obtain
an injective homomorphism $\C[v_i^{(n)}] \arr \C[\N \backslash
N_+/A_+]$. The fact that it is an isomorphism follows in the same way
as in the proof of \thmref{iso}, from the computation of
characters. Clearly,
\begin{equation}    \label{char}
\on{ch} \C[v_i^{(n)}] = \prod_{i=1}^\el \prod_{n_i>d_i}
(1-q^{n_i})^{-1}.
\end{equation}
On the other hand, since $\N$ acts freely on $N_+/A_+$,
$$
\on{ch} \C[\N \backslash N_+/A_+] = \on{ch} \C[N_+/A_+] (\on{ch}
\C[\N])^{-1},
$$
and
$$
\on{ch} \C[\N] = \prod_{i=1}^\el \prod_{n_i=1}^{d_i} (1-q^{n_i})^{-1},
$$
by \cite{Ko} (here we use $-\De$ as the gradation operator). Hence
$\on{ch} \C[\N \backslash N_+/A_+]$ is also given by formula
\eqref{char}, and the theorem is proved.
\end{proof}

\begin{rem}
One can prove that $\C[\N \backslash N_+/A_+]$ is isomorphic to a ring
of differential polynomials without using formula \eqref{vik} for the
functions $v_i$. However, the proof given above is shorter and more
explicit.\qed
\end{rem}

\subsection{KdV hierachy}
\thmref{vi} means that just like $N_+/A_+$, the double quotient $\N
\backslash N_+/A_+$ can be identified with the space of $\infty$--jets
of an $\el$--tuple of functions. We call this space the KdV jet space
and denote it by $\VV$. The vector fields $p_{-n}^R$ still act on $\N
\backslash N_+/A_+$ and hence give us an infinite set of commuting
evolutionary derivations on $\C[\VV]$. This is, by definition, the
{\em KdV hierarchy} associated to $\g$.

The natural projection $N_+/A_+ \arr \N \backslash N_+/A_+$ gives us a
map $\UU \arr \VV$, which amounts to expressing each $v_i$ as a
differential polynomial in $u_j$'s. This map is called the generalized
{\em Miura transformation}. It can be thought of as a change of
variables transforming the mKdV hierarchy into the KdV hierarchy.

\subsection{Drinfeld-Sokolov reduction}

In this section we give another definition of the generalized KdV hierarchy
following Drinfeld and Sokolov. We will show in the next section that the
two definitions are equivalent.

Denote by $\QQ$ the space of $\infty$--jets of functions ${\mathbf q}:
{\mathcal D}: \arr \bo_+$, where $\bo_+$ is the finite-dimensional
Borel subalgebra $\h \oplus \ovl{\n}_+$ (recall that $\ovl{\n}_+$ is
generated by $e_i, i=1,\ldots,\el$). If we choose a basis $\{
\omega^\vee_i, i=1,\ldots,\el \} \cup \{ e_\al, \al \in \Delta_+ \}$
of $\bo_+$, then we obtain the corresponding coordinates $q_\al^{(n)},
q_i^{(n)}, n\geq 0$ on $\QQ$. Thus, $\C[\QQ]$ is a ring of
differential polynomials equipped with an action of $\pa_z$.

Let $C(\g)$ be vector the space of operators of the form
\begin{equation}    \label{kdvL}
\pa_z + p_{-1} + {\mathbf q},
\end{equation}
where
$$
{\mathbf q} = \sum_{i=1}^\el \omega^\vee_i \otimes q_i + \sum_{\al
\in \Delta_+} e_\al \otimes q_\al \in \bo_+ \otimes \C[\QQ].
$$
The group $\N[\QQ]$ acts naturally on $C(\g)$:
\begin{equation}    \label{dei}
x \cdot (\pa_z + p_{-1} + {\mathbf q}) = x (\pa_z + p_{-1} + {\mathbf
q}) x^{-1} = \pa_z + [x,p_{-1} + {\mathbf q}] - \pa_z x.
\end{equation}
Note that $$p_{-1} = \pp + e_{\on{max}} \otimes t^{-1},$$ where
\begin{equation}    \label{ovlp}
\ovl{p}_{-1} = \sum_{j=1}^\el \frac{(\al_i,\al_i)}{2} f_j \in
\ovl{\n}_-
\end{equation}
and $e_{\on{max}}$ is a generator of the one-dimensional subspace of
$\g$ corresponding to the maximal root.

We have a direct sum decomposition $\bo_+ = \oplus_{i\geq 0}
\bo_{+,i}$ with respect to the principal gradation. The operator
$\on{ad} \pp$ acts from $\bo_{+,i+1}$ to $\bo_{+,i}$ injectively for
all $i>1$. Hence we can choose for each $j>0$ a vector subspace $S_j
\subset \bo_{+,j}$, such that $\bo_{+,j} = [\pp,\bo_{+,{j+1}}] \oplus
S_j$. Note that $S_j \neq 0$ if and only if $j$ is an exponent of
$\g$, and in that case $\dim S_j$ is the multiplicity of the exponent
$j$. In particular, $S_0=0$. Let $S = \oplus_{j\in E} S_j \subset
\nn_+$, where $E$ is the set of exponents of $\g$. Then, by
construction, $S$ is transversal to the image of the operator $\on{ad}
\ovl{p}_{-1}$ in $\ovl{\n}_+$.

\begin{prop}[\cite{DS}]    \label{free}
The action of $\N[\QQ]$ on $C(\g)$ is free. Furthermore, each
$\N[\QQ]$--orbit contains a unique operator \eqref{kdvL} satisfying
the condition that ${\mathbf q} \in S[\QQ]$.
\end{prop}

\begin{proof}
Denote by $C'(\g)$ the space of operators of the
form
$$
\pa_z + \pp + {\mathbf q}, \quad \quad {\mathbf q} \in \bo_+ \otimes \C[\QQ].
$$
The group $\N[\QQ]$ acts on $C'(\g)$ by the formula analogous to
\eqref{dei}. We have an isomorphism $C(\g) \arr C'(\g)$, which sends
$\pa_z + p_{-1} + {\mathbf q}$ to $\pa_z + \pp + {\mathbf q}$. Since
$[x,e_{\on{max}}] = 0, \forall x \in \nn_+$, this isomorphism commutes
with the action of $\N[\QQ]$. Thus, we can study the action of $\N[\QQ]$ on
$C'(\g)$ instead of $C(\g)$.

We claim that each element of $C(\g)$ can be uniquely represented in
the form
\begin{equation}    \label{gauge}
\pa_z + \pp + {\mathbf q} = \exp \left( \on{ad} U \right) \cdot \left(
\pa_z + \pp + {\mathbf q}^0 \right),
\end{equation}
where $U \in \nn_+[\QQ]$ and ${\mathbf q}^0 \in S[\QQ]$. Decompose with
respect to the principal gradation: $U=\sum_{j\geq 0} U_j$, ${\mathbf
q} = \sum_{j\geq 0} {\mathbf q}_j$, ${\mathbf q}^0 = \sum_{j\geq 0}
{\mathbf q}^0_j$. Equating the homogeneous components of degree $j$ in
both sides of \eqref{gauge}, we obtain that ${\mathbf q}^0_i +
[U_{i+1},\pp]$ is expressed in terms of ${\mathbf q}_i,{\mathbf
q}^0_1,\ldots,{\mathbf q}^0_{i-1},U_1,\ldots,U_i$. The direct sum
decomposition $\bo_{+,i} = [\pp,\bo_{+,{i+1}}] \oplus S_i$ then allows
us to determine uniquely ${\mathbf q}^0_i$ and $U_{i+1}$. Hence $U$
and ${\mathbf q}^0$ satisfying equation \eqref{gauge} can be found
uniquely by induction, and lemma follows.
\end{proof}

{}From the analytic point of view, \propref{free} means that every
first order differential operator \eqref{kdvL}, where ${\mb q}: \R
\arr \bo_+$ is a smooth function, can be brought to the form
\eqref{kdvL}, where ${\mb q}$ takes values in $S \subset \nn_+$, by
gauge transformation with a function $x(z): \R \arr \N$. Moreover,
$x(z)$ depends on the entries of ${\mb q}$ only through their
derivatives at $z$ (i.e., it is local).

\begin{rem}
The statement of the lemma remains true if we replace the ring of
differential polynomials $\C[\QQ]$ by any differential ring $R$. For
example, we can take $R = \C[[z]]$. Then the quotient of $C(\g)$ by
the action of $\N((z))$ is what Beilinson and Drinfeld call the space
of {\em opers} on the formal disc. This space can be defined
intrinsically without choosing a particular uniformizing parameter $z$
(using the notion of connection on the formal disc). In this form it
has been generalized by Beilinson and Drinfeld to the situation where
the formal disc is replaced by any algebraic curve. For instance, in
the case of $\g=\sw_2$, the notion of oper coincides with that of
projective connection.\qed
\end{rem}

It is easy to see from the proof of \propref{conj} that its statement
remains true if we replace the operator $\pa_z + p_{-1} + {\mathbf
u}(z)$ by operator \eqref{kdvL}. Using this fact, Drinfeld and Sokolov
construct in \cite{DS} the zero-curvature equations for the operator
\eqref{kdvL} in the same way as for the mKdV hierarchy using formulas
\eqref{ds}. Drinfeld and Sokolov show that these equations preserve
the corresponding $\N[\QQ]$--orbits (see \cite{DS}, Sect.~6.2). Thus,
they obtain a system of compatible evolutionary equations on
$\N[\QQ]$--orbits in $C(\g)$, which they call the generalized KdV
hierarchy corresponding to $\G$.

These equations give rise to evolutionary derivations acting on the
ring of differential polynomials $\C[s_i^{(n)}]_{i=1,\ldots,\el;n\geq
0}$. Indeed, the space $S$ is $\el$--dimensional. Let us choose
homogeneous coordinates $s_1,\ldots,s_\el$ of $S$. Then according to
\propref{free}, the KdV equations can be written as partial
differential equations on $s_i$'s. It is easy to see that the first of
them is just $\pa_z$ itself. Hence others give rise to evolutionary
derivations of $\C[s_i^{(n)}]_{i=1,\ldots,\el;n\geq 0}$.

\subsection{Equivalence of two definitions}

Recall that $\QQ$ is the space of $\infty$--jets of ${\mathbf q}$ with
coordinates $q_\al^{(n)}, q_i^{(n)}, n\geq 0$. Denote by $\Ss$ the
space of $\infty$--jets of ${\mathbf q}^0$ with coordinates
$s_i^{(n)}, i=1,\ldots,\el, n\geq 0$.

Note that we have a natural embedding $\imath: \UU \arr \QQ$, which
sends $u_i^{(n)}$ to $q_i^{(n)}$. Let $\mu: \UU \arr \Ss$ be the
composition of the embedding $\imath$ and the projection $\QQ \arr
\Ss$.

\propref{conj}, suitably modified for operators of the form
\eqref{kdvL}, gives us a map $\nu: \QQ \arr N_+/A_+$ and hence a map
$\wt{\nu}: \Ss \arr \N \backslash N_+/A_+$. According to
\thmref{ident}, the composition of the embedding $\imath$ with the map
$\nu$ coincides with the isomorphism $\UU \simeq N_+/A_+$ constructed
in \thmref{iso}. Hence the maps $\nu$ and $\wt{\nu}$ are surjective.

\begin{prop}    \label{est}
The map $\wt{\nu}: \Ss \arr \N \backslash N_+/A_+$ is an isomorphism.
\end{prop}

\begin{proof}
By construction, the homomorphism $\wt{\nu}$ is surjective and
homogeneous with respect to the natural $\Z$--gradations. Hence, it
suffices to show that the characters of the spaces $\C[\Ss]$ and
$\C[\N \backslash N_+/A_+]$ coincide. Since, by construction, $\deg
s_i = d_i + 1$, $$\on{ch} \C[\Ss] = \prod_{i=1}^\el \prod_{n_i \geq
d_i+1} (1-q^{n_i})^{-1}.$$ But in the proof of \thmref{vi} we showed
that that $\on{ch} \C[\N \backslash N_+/A_+]$ is given by the same
formula.
\end{proof}

Thus, we have shown that there is an isomorphism of rings
$$\C[v_i^{(n)}]_{i=1,\ldots,\el;n\geq 0} \arr
\C[s_i^{(n)}]_{i=1,\ldots,\el;n\geq 0},$$ which preserves the
$\Z$--gradation and the action of $\pa_z$, and such that it sends the
derivations of generalized KdV hierarchy defined in Sect.~4.3 to the
derivations defined in Sect.~4.4 following Drinfeld-Sokolov
\cite{DS}. It also follows that the map $\mu$ defined above as the
composition of the embedding $\UU \arr \QQ$ and the projection $\QQ
\arr \N \backslash \QQ \simeq \Ss$ coincides with the projection $\UU
\simeq N_+/A_+ \arr \N \backslash N_+/A_+ \simeq \Ss$. This is the
generalized Miura transformation.

\subsection{Example of $\sw_2$}

The space $C'(\sw_2)$ consists of operators of the form
$$
\pa_z + \begin{pmatrix} a & b \\ 1 & -a \end{pmatrix}.
$$
The group $$\N = \begin{pmatrix} 1 & x \\ 0 & 1 \end{pmatrix}$$ acts on
$C'(\sw_2)$ by gauge transformations \eqref{dei}. We write canonical
representatives of $\N$--orbits in the form
$$
\begin{pmatrix} 0 & s \\ 1 & 0 \end{pmatrix}.
$$
We have:
$$
\begin{pmatrix} 1 & - \frac{u}{2} \\ 0 & 1 \end{pmatrix} \left( \pa_z +
\begin{pmatrix} \frac{u}{2} & 0 \\ 1 & - \frac{u}{2} \end{pmatrix} \right)
\begin{pmatrix} 1 & - \frac{u}{2} \\ 0 & 1 \end{pmatrix}^{-1} = 
\begin{pmatrix} 0 & s \\ 1 & 0 \end{pmatrix},
$$
where
\begin{equation}    \label{v}
s = \frac{1}{4} u^2 + \frac{1}{2} \pa_z u.
\end{equation}
This formula defines the Miura transformation $\C[s^{(n)}] \arr
\C[u_i^{(n)}]$. If we apply this change of variables \eqref{v} to the mKdV
equation \eqref{mkdv1} we obtain the equation
\begin{equation}    \label{kdv1}
\pa_3 s = \frac{3}{2} s \pa_z s - \frac{1}{4} \pa_z^3 s,
\end{equation}
which, as expected, closes on $s$ and its derivatives. This equation
becomes the KdV equation \eqref{kdv} after a slight redefinition of
variables: $s \arr -v$, $\tau_3 \arr - 4 \tau$. Thus, the results of
this lecture prove the existence of infinitely many higher KdV flows,
i.e., evolutionary derivations on $\C[v^{(n)}]$, which commute with
the KdV derivation defined by equation \eqref{kdv}.

The fact that the KdV and mKdV equations are connected by a change of
variables \eqref{kdv1} was discovered by R.~Miura, who also realized that
it can be rewritten as
$$
\pa_z^2 - s = \left( \pa_z - \frac{u}{2} \right) \left( \pa_z + \frac{u}{2}
\right).
$$
It appears that it was this observation that triggered the fascinating idea
that the KdV equation should be considered as a flow on the space of
Sturm-Liouville operators $\pa^2_z - s(z)$ \cite{GGKM}, which led to the
concept of Inverse Scattering Method and the modern view of the theory of
solitons (see \cite{Ne}).

\subsection{Explicit formulas for the action of $\n_+$}    \label{expl}

In this section we obtain explicit formulas for the action of the
generators $e_i, i=0,\ldots,\el$, of $\n_+$ on $\C[\UU]$. This formulas can
be used to find the KdV variables $v_i \in \C[\UU]$, and we will also need
them in the next lecture when we study Toda field theories.

In order to find these formulas, we need a geometric construction of
modules contragradient to the Verma modules over $\G$ and
homomorphisms between them. This construction is an affine analogue of
Kostant's construction \cite{Kos} in the case of simple Lie algebras.

For $\la \in \h^*$, denote by $\C_\la$ the one-dimensional
representation of ${\mathfrak b}_+$, on which $\h \subset {\mathfrak
b}_+$ acts according to its character $\la$, and $\n_+ \subset
{\mathfrak b}_+$ acts trivially. Let $M_\la$ be the Verma module over
$\G$ of highest weight $\la$: $$M_\la = U(\G) \otimes_{U({\mathfrak
b}_+)} \C_\la.$$ we denote by $v_\la$ the highest weight vector of
$M_\la^*$, $1 \otimes 1$.

Denote by $\lan\cdot,\cdot\ran$ the pairing $M_\la^* \times M_\la
\arr \C$. Let $\omega$ be the Cartan anti-involution on $\G$, which maps
generators $e_0,\ldots,e_\el$ to $f_0,\ldots,f_\el$ and vice versa and
preserves $\h$ \cite{Kac}. It extends to an anti-involution of $U(\G)$. Let
$M^*_\la$ be the module contragradient to $M_\la$. As a linear space,
$M_\la^*$ is the restricted dual of $M_\la$. The action of $x \in \G$ on $y
\in M_\la^*$ is defined as follows: $$\lan x\cdot y,z \ran = \lan
y,\omega(x)\cdot z \ran, \quad \quad z \in M_\la.$$

The module $M_\la^*$ can be realized in the space of regular functions on
$N_+$ (see \cite{FF:kdv}, Sect.~4). Indeed, the $\n_+$--module $\C[N_+]$
(with respect to the right action) is dual to a free module with one
generator, and so is each $M_\la^*$. Hence we can identify $M_\la^*$ and
$\C[N_+]$ as $\n_+$--modules for any $\la$. It is easy to see that for
$\la=0$, $M_0^*$ is isomorphic to $\C[N_+]$ on which the action of $\G$ is
defined via the left infinitesimal action of $\G$ on $N_+$ by vector
fields, described in Lecture 1. For general $\la$, the action of $\G$ is
given by first order differential operators: for $a \in \G$ this
differential operator is equal to $a^R + f_\la(a)$, where $f_\la(a) \in
\C[N_+]$. The function $f_\la$ is the image of $a \cdot v_\la$ under the
isomorphism $M_\la^* \simeq \C[N_+]$. Hence if $a$ is homogeneous, then
$f_\la(a)$ is also homogeneous of the same weight.

Here is another, homological, point of view on $f_\la(a)$. As a
$\G$--module, $\C[N_+] = M_0^*$ is coinduced from the trivial
representation of ${\mathfrak b}_-$. By Shapiro's lemma (cf., e.g.,
\cite{fuchs}, Sect.~1.5.4), $H^1(\G,\C[N_+]) \simeq H^1({\mathfrak
b}_-,\C) \simeq ({\mathfrak b}_-/[{\mathfrak b}_-,{\mathfrak b}_-])^*
= \h^*$. We see that all elements of $H^1(\G,\C[N_+])$ have weight
$0$. On the other hand, functions on $N_+$ can only have negative or
$0$ weights and the only functions, which have weight $0$ are
constants, which are invariant with respect to the action of
$\G$. Therefore the coboundary of any element of the $0$th group of
the complex, $\C[N_+]$, has a non-zero weight. Hence any cohomology
class from $H^1(\G,\C[N_+])$ canonically defines a one-cocycle $f$,
i.e., a map $\G \arr \C[N_+]$. Thus, having identified
$H^1(\G,\C[N_+])$ with $\h^*$, we can assign to each $\la \in \h^*$
and each $a \in \G$, a function on $N_+$ -- this is our $f_\la(a)$.

The following two results will enable us to compute the action of $e_i^L$.

\begin{prop}[\cite{Kos}, Theorem 2.2]    \label{kost}
Consider $\la \in \h^*$ as an element of $\h$, using the invariant inner
product. Then
\begin{equation}    \label{fla}
f_\la(a)(x) = (\la,x a x^{-1}).
\end{equation}
\end{prop}

\begin{prop}[\cite{FF:kdv}, Prop.~3]    \label{P}
For any $a \in \G$ we have:
\begin{equation}    \label{key}
[e_i^L,a^R] = - f_{\al_i}(a) e_i^L, \quad \quad i=0,\ldots,\el.
\end{equation}
\end{prop}

Substituting $\la=\al_i$ and $a=p_{-1}$ in formula \eqref{fla}, we obtain:
$$
f_{\al_i}(p_{-1})(x) = (\al_i,x p_{-1} x^{-1}) = u_i.
$$
Hence formula \eqref{key} gives:
\begin{equation}    \label{ui}
[e_i^L,p_{-1}^R] = - u_i e_i^L.
\end{equation}
Let us write
$$
e_i^L = \sum_{1\leq j\leq l, n\geq 0} C_{i,j}^{(m)} \fp{j}{m},
$$
where $C_{i,j}^{(m)} \in \C[u_i^{(n)}]_{1\leq i\leq \el;n\geq 0}$.  Formula
\eqref{ui} now gives us recurrence relations for the coefficients of
$\pa/\pa u_j^{(m-1)}$ in the vector field $e_i$: $$C_{i,j}^{(m)} =
-u_i C_{i,j}^{(m-1)} + \pa_z C_{i,j}^{(m-1)},$$ where we have
identified $p_{-1}^R$ with $\pa_z$ and used the formula
$$\left[ \frac{\pa}{\pa u_i^{(n)}},\pa_z \right] = \frac{\pa}{\pa
u_i^{(n-1)}},$$ if $n>0$. We also have, according to formula \eqref{neg},
$$e_i \cdot u_j = - (\al_i,\al_j).$$ This gives the initial condition for
our recurrence relation. Combining them, we obtain the following formula:
\begin{equation}    \label{eil}
e_i^L = - \sum_{n\geq 0} B_i^{(n)} \pa_i^{(n)},
\end{equation}
where
\begin{equation}    \label{pain}
\pa_i^{(n)} = \sum_{1\leq j\leq l} (\al_i,\al_j) \frac{\pa}{\pa
u_j^{(n)}},
\end{equation}
and $B_i^{(n)}$'s are polynomials in $u_i^{(m)}$'s, which satisfy the
recurrence relation:
\begin{equation}    \label{rec}
B_i^{(n)} = -u_i B_i^{(n-1)} + \pa_z B_i^{(n-1)},
\end{equation}
with the initial condition $B_i^{(0)} = 1$. 

Formula \eqref{eil} is valid for all $i=0,\ldots,\el$. Using this formula,
we can find explicitly the KdV variables $v_i, i=1,\ldots,\el$, following
the proof of \thmref{vi}. Consider, for example, the case $\g=\sw_2$. In
this case $\C[\UU] = \C[u^{(n)}]$, and $\C[\VV] = \C[v^{(n)}]$ is the
subspace of $\C[\UU]$, which consists of differential polynomials
annihilated by
$$
e^L = - \sum_{n\geq 0} B^{(n)} \frac{\pa}{\pa u^{(n)}}.
$$
We find from formula \eqref{rec}: $B^{(0)} = 1, B^{(1)} = -u$, etc.

The KdV variable $v$ lies in the degree $2$ component of $\C[\UU]$, i.e.,
$\on{span} \{ u^2,\pa_z u \}$, and we have: $e^L \cdot u^2 = -2u, e^L \cdot
\pa_z u = u$. Hence
\begin{equation}    \label{v1}
v = \frac{1}{4} u^2 + \frac{1}{2} \pa_z u
\end{equation}
is annihilated by $e^L$ and can be taken as a KdV variable. As expected,
this agrees with formula \eqref{v} obtained by means of the
Drinfeld-Sokolov reduction.

\setcounter{section}{5}
\setcounter{subsection}{0}
\setcounter{equation}{0}
\setcounter{thm}{0}
\setcounter{prop}{0}
\setcounter{lem}{0}
\setcounter{rem}{0}

\section*{Lecture 5}

In this lecture we consider the Toda equations and their local integrals of
motion. The Toda equation associated to the simple Lie algebra $\g$
(respectively, affine algebra) reads
\begin{equation}   \label{toda}
\pa_\tau \pa_z \phi_i(z,t) = \sum_{j \in J} (\al_i,\al_j) e^{-\phi_j(z,t)},
\quad i = 1,\ldots,\el,
\end{equation}
where $J$ is the set of vertices of the Dynkin diagram of $\g$
(respectively, $\G$). Each $\phi_i(z,t), i=1,\ldots,\el$, is a family of
functions in $z$, depending on the time variable $\tau$, and $$\phi_0(z,t)
= - \frac{1}{a_0} \sum_{i=1}^\el a_i \phi_i(t,\tau).$$ If we set formally
$u_i(t) = \pa_z \phi_i(t), i=1,\ldots,\el$, then equations \eqref{toda} can
be rewritten as an equation on ${\mathbf u}(z)$:
\begin{equation}   \label{toda1}
\pa_\tau u_i = \sum_{j \in J} (\al_i,\al_j) e^{-\phi_j}, \quad i =
1,\ldots,\el,
\end{equation}
where $\phi_j$ is now understood as the anti-derivative of
$u_j$.

Certainly, equation \eqref{toda1} is not in the form
\eqref{higherkdv}, since the terms $e^{-\phi_j}$ on the right hand
side are not differential polynomials in $u_i$'s, so our differential
algebra formalism can not be applied. However, we can modify our
formalism to suit the Toda equations. Moreover, we will see that these
equations are closely connected to the derivations $e_i^L$ studied in
the previous lecture. After that, we will introduce the notion of
local integral of motion of the Toda equation. We will describe these
local integrals of motion, and will show that the corresponding
hamiltonian flows are the mKdV flows.

\subsection{The differential algebra formalism of Toda equation}

Let $\La$ be the $(\el+1)$ dimensional lattice spanned by
$\al_0,\ldots,\al_\el$. For any element $\la = \sum_{0\leq i\leq l}
\la_i \al_i$ of $\La$, define the linear space $\pi_\la =
\C[u_i^{(n)}] \otimes e^{\laa}$, where $\ovl{\la} = \sum_{0\leq i\leq
l} \la_i \phi_i$, equipped with an action of $\pa_z$ by the formula
\begin{equation}    \label{pala}
\pa_z \cdot (P \otimes e^{\laa}) = (\pa_z P) \otimes e^{\laa} + \left(
\sum_{0\leq i\leq l} \la_i u_i^{(0)} P \right) \otimes e^{\laa},
\end{equation}
where we put $$u_0^{(n)} = -\frac{1}{a_0} \sum_{1\leq i\leq l} a_i
u_i^{(n)}.$$ This formula means that we set $\pa_z \phi_i = u_i$. In
particular, $\pi_0 = \C[u_i^{(n)}]$.

We introduce a $\Z$--gradation on $\pi_\la$ by setting $\deg P \otimes
e^{\laa} = (\rho^\vee,\la) + \deg P$, where the gradation on
$\C[u_i^{(n)}]$ is the one previously defined (recall that $\rho^\vee$
satisfies $(\rho^\vee,\al_i) = 1, i=0,\ldots,\el$).

An evolutionary operator $X: \pi_\la \arr \pi_\mu$ is, by definition, a
linear operator $\pi_\la \arr \pi_\mu$, which commutes with the action of
$\pa_z$. Defining an evolutionary operator $X$ acting from $\pi_0$ to
$\oplus_{j \in S} \pi_{\la_j}$ is the same as defining $X \cdot u_i$ as an
element of $\oplus_{j \in S} \pi_{\la_j}$ for all $j \in S$. Hence
given an equation of the form
$$\pa_\tau u_i = \sum_{j \in S} P_{ij},
$$
where $P_j \in \pi_{\la_j}$, we obtain an evolutionary operator
$$
X\{P_{ij}\} = \sum_{j \in S} X_j,
$$
where
$$
X_j = \sum_{n\geq 0} \sum_{i=1}^\el (\pa_z^n \cdot P_{ij}) \fp{i}{n}.
$$

Toda equation \eqref{toda1} gives rise to the evolutionary operator
\begin{equation}
{\mathcal H}: \pi_0 \arr \oplus_{j \in J} \pi_{-\al_j},
\end{equation}
which is the sum of the terms
\begin{equation}    \label{borel}
\q_i = \sum_{n\geq 0} (\pa_z^n e^{-\phi_i}) \pa_i^{(n)},
\end{equation}
where $\pa_i^{(n)}$ is given by formula \eqref{pain}. Here $\q_i$ is
an evolutionary operator acting from $\pi_0$ to
$\pi_{-\al_i}$. Clearly,
$$\pa_z^n e^{-\phi_i} = b_i^{(n)} e^{-\phi_i},$$ where $b_i^{(n)}$'s are
differential polynomials in $u_j^{(n)}$'s, satisfying the recurrence
relation
$$
b_i^{(n)} = -u_i b_i^{(n-1)} + \pa_z b_i^{(n-1)}.
$$
This recurrence relation coincides with the recurrence relation \eqref{rec}
for the differential polynomials $B_i^{(n)}$ appearing in formula
\eqref{eil} for $e_i^L$. Hence we obtain the following result.

\begin{prop}
Let $T_j$ be the operator of multiplication by $e^{-\phi_j}$ acting
from $\pi_0$ to $\pi_{-\al_i}$. The evolutionary operator ${\mathcal
H}$ associated to the Toda equation equals
$$
{\mathcal H} = - \sum_{j \in J} T_j e_j^L.
$$
\end{prop}

Thus, while the mKdV flows come from the right action of $\ab_+$ on
$N_+/A_+$, the Toda flow comes from the left action of the generators
of $\n_+$ (or $\nn_+$) on $N_+/A_+$.

Now we want to introduce the notion of local integral of motion of
the Toda equation. In order to do that, we first review the general
hamiltonian formalism of soliton equations.

\subsection{Hamiltonian formalism}

In this section we briefly discuss the hamiltonian structure on the
space of local functionals. This is a special case of the formalism of
generalized hamiltonian structures and of generalized hamiltonian
operators developed by Gelfand, Dickey and Dorfman
\cite{gd1,gd2,gdorf} (see \cite{FF:laws}, Sect.~2; \cite{FF:kdv},
Sect.~2 for review).

First we introduce the concept of {\em local functional}. Consider the
space $F(\h)$ of functions on the circle with values in the Cartan
subalgebra $\h$ of $\g$, ${\mathbf u}(z) =
(u_1(z),\ldots,u_\el(z))$. Then each differential polynomial $P \in
\C[u_i^{(n)}]$ gives rise to a functional on $F(\h)$ that
sends ${\mathbf u}(z)$ to
$$
\int P(u_i^{(n)}) \; dz.
$$
Such functional are called local functionals. It is easy to see that
this functional is $0$ if and only if $P$ is a sum of a total
derivative (i.e., $\pa_z Q$ for some $Q \in \C[u_i^{(n)}]$) and a
constant. Hence we formally define the space of local functionals
$\F_0$ as the quotient of $\pi_0=\C[u_i^{(n)}]$ by the subspace
$\on{Im} \pa_z \oplus \C$ spanned by total derivatives and
constants. We introduce a $\Z$-gradation on $\F_0$ by subtracting $1$
from the gradation induced from $\pi_0$, and denote by $\int$ the
projection $\pi_0 \arr \F_0$.

From now on, to simplify notation, we shall write $\pa$ for $\pa_z$.

Introduce variational derivatives: $$\delta_i P = \sum_{n\geq 0}
\sum_{1\leq j\leq l} (\al_i,\al_j) (-\pa)^n \frac{\pa P}{\pa
u_j^{(n)}}, \quad \quad i=1,\ldots,\el.$$ Now for each $P \in \pi_0$,
define an evolutionary derivation $\xi_P$ of $\pi_0$ by the formula.
\begin{equation}    \label{vf}
\xi_P = \sum_{1\leq i\leq l, n\geq 0} ( \pa^{n+1} \cdot \delta_i P)
\fp{i}{n},
\end{equation}
Since $\delta_i P = 0, i=1,\ldots,\el, \forall P \in \on{Im} \pa_z
\oplus \C$, we see that $\xi_P$ depends only on the image of $P$ in
$\F_0$, $\int P dz$. So sometimes we shall write $\xi_{\int P dz}$
instead of $\xi_P$. We can view $\xi: P \arr \xi_P$ as a linear map
$\F_0 \arr \on{Der}^{\pa_z}$, where $\on{Der}^{\pa_z}$ is the space of
evolutionary derivations on $\pi_0$ (note that it is a Lie algebra).

For instance, we find:
\begin{equation}    \label{pa}
\pa = \sum_{1\leq i\leq l, n\geq 0} u_i^{(n+1)}
\frac{\pa}{\pa u_i^{(n)}} = \xi_P, \quad \quad P = \frac{1}{2} \sum_{1\leq
i\leq l} u_i u^i,
\end{equation}
where $u^i, i=1\ldots,\el,$ are dual to $u_i, i=1\ldots,\el,$ with
respect to the inner product defined by $(\cdot,\cdot)$.

Now we define a Poisson bracket on $\F_0$ by the formula
\begin{equation}    \label{pb}
\{\int P dz,\int R dz\} = \int (\xi_P \cdot R) \; dz.
\end{equation}

\begin{prop}[\cite{gd1,gd2,gdorf}]
Formula \eqref{pb} defines a Lie algebra sructure on
$\F_0$. Furthermore, the map $\xi: \F_0 \arr \on{Der}^\pa$ is a
homomorphism of Lie algebras, so that
\begin{equation}    \label{com}
\xi_{\{ P , R\}} = [\xi_P , \xi_R], \quad \quad \forall P,R \in \pi_0.
\end{equation}
\end{prop}

\begin{rem}
Using the isomorphism $\pi_0 \simeq \C[N_+/A_+]$, B.~Enriquez and the
author have expressed the formula for Poisson bracket \eqref{pb} in
terms of the corresponding unipotent cosets \cite{EF:vpa}.\qed
\end{rem}

Now we extend the map $\xi$ to incorporate the spaces $\pi_\la$, following
Kuperschmidt and Wilson \cite{kw,W}.

Let $\F_\la$ be the quotient of $\pi_\la$ by the subspace of total
derivatives and $\int$ be the projection $\pi_\la \arr \F_\la, P \arr
\int P dz$. We define a $\Z$--gradation on $\F_\la$ by subtracting $1$
from the gradation induced from $\pi_\la$.

For any $P \in \F_0$ the derivation $\xi_P: \pi_0 \arr \pi_0$ can be
extended to a linear operator on $\oplus_{\la \in \La} \pi_\la$ by the
formula $$\xi_P = \sum_{1\leq i\leq l, n\geq 0} ( \pa^{n+1} \cdot \delta_i
P) \frac{\pa}{\pa u_i^{(n)}} + \sum_{1\leq i\leq l} \delta_i P
\frac{\pa}{\pa \phi_i},$$ where $\pa/\pa
\phi_i \cdot (S e^{\laa}) = \la_i S e^{\laa}$. This defines a structure
of $\F_0$-module on $\pi_\la$.

For each $P \in \pi_0$ the operator $\xi_P$ commutes with the action
of derivative. Hence we obtain the structure of an $\F_0$-module on
$\F_\la$, i.e., a map $\{ \cdot,\cdot \}: \F_0 \times \F_\la \arr
\F_\la$: $$\{ \int P dz,\int R dz \} = \int (\xi_{P} \cdot R) \; dz.$$

Similarly, any element $R \in \pi_\la$ defines a linear operator $\xi_R$,
acting from $\pi_0$ to $\pi_\la$ and commuting with $\pa$:
\begin{equation}    \label{oper}
\xi_{S e^{\laa}} = \sum_{1\leq i\leq l, n\geq 0} \pa^n \left( \pa
(\delta_i S \cdot e^{\laa}) -  S \, \frac{\pa e^{\laa}}{\pa
\phi_i} \right) \fp{i}{n}.
\end{equation}
The operator $\xi_R$ depends only on the image of $R$ in
$\F_\la$. Therefore it gives rise to a map $\{ \cdot,\cdot \}: \F_\la
\times \F_0 \arr \F_\la$. We have for any $P \in \F_0, R \in \F_\la$:
$$\int (\xi_R \cdot P) \; dz = - \int (\xi_P \cdot R) \; dz.$$
Therefore our bracket $\{ \cdot,\cdot \}$ is antisymmetric.

Note that, by construction, the operators $\xi_R$ are
evolutionary. Furthermore, formula \eqref{com} holds for any $P \in
\F_0, R \in \oplus_{\la \in \La} \F_\la$.

\subsection{Local integrals of motion}

We obtain from the definition of $\xi$:
$$
\xi_{e^{-\phi_j}} = \sum_{n\geq 0} (\pa^n e^{-\phi_i}) \pa_i^{(n)},
$$
which coincides with formula \eqref{borel} for $\q_i$. Hence
$$
\q_j = \left\{ \int e^{-\phi_j} dz, \cdot \right\},
$$
is a hamiltonian operator, and we obtain:

\begin{prop}
The evolutionary operator ${\mathcal H}$ defined by the Toda equation is
hamiltonian:
\begin{equation}    \label{mch}
{\mathcal H} = \sum_{j \in J} \left\{ \int e^{-\phi_j} dz, \cdot \right\}.
\end{equation}
\end{prop}

Consider the corresponding operator
$$
\{ H,\cdot \}: \F_0 \arr \oplus_{0\leq i\leq l}
\F_{-\al_i},
$$
where
$$
H = \sum_{j \in J} \int e^{-\phi_j} dz
$$
can be viewed as the hamiltonian of the Toda equation \eqref{toda1}.

We define {\em local integrals of motion} of the (finite of affine)
Toda field theory as local functionals, which lie in the kernel of the
operator ${\mathcal H}$. Thus, local integrals of motion are
quantities that are conserved with respect to the Toda equation.

Denote the space of all local integrals of motion associated to $\g$
by $I(\g)$ and the space of local integrals of motion associated to
$\G$ by $I(\G)$. Denote by $\Q_i$ the operator $\F_0 \arr \F_{-\al_i}$
induced by $\q_i$. Then,
$$
I(\g) = \bigcap_{i=1}^\el \on{Ker}_{\F_0} \Q_i, \quad \quad I(\G) =
\bigcap_{i=0}^\el \on{Ker}_{\F_0} \Q_i.
$$
Thus, we have a sequence of embeddings: $I(\G) \subset I(\g) \subset
\F_0$. From the fact that $\q_i$'s are hamiltonian operators, it
follows that both $I(\G)$ and $I(\g)$ are Poisson subalgebras of
$\F_0$.

Recall that in the previous lecture we defined an embedding
$\C[v_i^{(n)}] \arr \C[u_i^{(n)}]$. Let $\W(\g)$ be the quotient of
$\C[v_i^{(n)}]$ by the total derivatives and constants. We have an
embedding $\W(\g) \arr \F_0$.

The following results describe the spaces $I(\g)$ and $I(\G)$.

\begin{thm}[\cite{FF:laws}, (2.4.10)]    \label{walg}
{\em The space of local integrals of motion of the Toda equation associated
to a simple Lie algebra $\g$ equals $\W(\g)$.}
\end{thm}

\begin{thm}[\cite{FF:kdv}, Theorem 1]    \label{span}
{\em The space of local integrals of motion of the Toda equation
associated to an affine Kac-Moody algebra $\G$ is linearly spanned by
elements $I_m \in \F_0, m \in I$, where $\deg I_m = m$.}
\end{thm}

The proof of these results is based on the following idea
\cite{FF:laws}. We need to compute the $0$th cohomology of the two-step
complex
$$
\F_0 \arr \oplus_{j \in J} \F_{-\al_j},
$$
where the differential is given by the sum of $\Q_i, i\in S$. We will
extend this complex further to the right (this certainly does not
change its $0$th cohomology), and then we will compute the cohomology
of the resulting complex by using the Bernstein-Gelfand-Gelfand (BGG)
resolution.

\subsection{The complexes $F^\bu(\G)$ and $F^\bu(\g)$}    \label{bggres}

Recall \cite{bgg,rocha} that the dual of the BGG resolution of $\G$ is
a complex $B^\bu(\G) = \oplus_{j\geq 0} B^j(\G)$, where $B^j(\G) =
\oplus_{w \in W, l(w)=j} M^*_{w(\rho) - \rho}$. Here $M^*_\la$ is the
module contragradient to the Verma module of highest weight $\la$, $W$
is the corresponding affine Weyl group, and $l: W \arr \Z_+$ is the
length function. The differentials of the complex are described in the
Appendix. They commute with the action of $\G$. The $0$th cohomology
of $B^\bu(\G)$ is one-dimensional and all higher cohomologies of
$B^\bu(\G)$ vanish, so that $B^\bu(\G)$ is an injective resolution of
the trivial representation of $\n_+$.

By construction, the right action of the Lie algebra $\G$ on this
complex commutes with the differentials. Therefore we can take the
subcomplex of invariants $F^\bu(\G) = B^\bu(\G)^{\ab_+^R}$ with
respect to the action of the Lie algebra $\ab_+^R$. Since $M_\la^*
\simeq \C[N_+]$ as an $\n_+$--module (see Sect.~4.7), we can identify
$(M_\la^*)^{\ab_+^R} \simeq \C[N_+/A_+]$ with $\pi_\la = \C[u_i^{(n)}]
\otimes e^{\laa}$ using \thmref{iso}. Hence
$$F^0(\G) = \pi_0, \quad \quad F^1(\G) = \oplus_{i=0}^\el
\pi_{-\al_i}.$$ Moreover, according to Sect.~A.2 of the Appendix, the
differential $F^0(\G) \arr F^1(\G)$ equals ${\mathcal H}$ given by
formula \eqref{mch} with $J=\{ 0,\ldots,\el \}$.

By construction, the differential of the complex $F^\bu(\G)$ commutes
with the action of $\ab_-^R$, in particular, with the action of
$\pa=p_{-1}^R$. It is also clear that the space of constants $\C
\subset \pi_0 = F^0(\G)$ is a direct summand in $F^\bu(\G)$. Hence the
quotient $F^\bu(\G)/(\C \oplus \on{Im} \pa)$ is also a complex. We have,
by definition, $$\F^0(\G) = \F_0, \quad \quad \F^1(\G) = \oplus_{0\leq
i\leq \el} \F_{-\al_i},$$ and the differential $\bar{\delta}^0: \F^0(\G)
\arr \F^1(\G)$ is given by $$\bar{\delta}^0 = \sum_{0\leq i\leq \el}
\q_i.$$ Therefore the $0$th cohomology of the complex $\F^\bu(\G)$ is
isomorphic to the $I(\G)$. This cohomology will be computed in the
next section.

\subsection{Proof of \thmref{span}}

First we compute the cohomology of the complex $F^\bu(\G)$.

\begin{lem}[\cite{FF:kdv}, Lemma 1]    \label{action}
The cohomology of the complex $F^\bu(\G)$ is isomorphic to
$H^\bu(\n_+,\pi_0)$. The action of $\ab_-^R$ on the cohomology of
$F^\bu(\G)$ is trivial.
\end{lem}

Here $H^\bu(\n_+,\pi_0)$ denotes the Lie algebra cohomology of $\n_+$
with coefficients in $\pi_0 \simeq \C[N_+/A_+]$ (see \cite{fuchs} for
the definition). The $\n_+$--module $\C[N_+/A_+]$ is coinduced from
the trivial representation of $\ab_+$. Hence, by Shapiro's lemma (see
\cite{fuchs}, Sect.~1.5.4), $H^\bu(\n_+,\pi_0) \simeq
H^\bu(\ab_+,\C)$. Since $\ab_+$ is an abelian Lie algebra,
$H^\bu(\ab_+,\C)$ is isomorphic to the exterior algebra
$\bigwedge^\bu(\ab_+^*) \simeq \bigwedge^\bu(\ab_-)$.

\begin{lem}    \label{ide}
The $0$th cohomology of the complex $\F^\bu(\G)$ is isomorphic to the
$1$st cohomology of the complex $F^\bu(\G)$, and hence to $\ab_-$.
\end{lem}

This readily implies \thmref{span}. Indeed, with respect to the
$\Z$--gradation on $\F_\la$, the differentials of the complex are
homogeneous of degree $0$. Moreover, the corresponding $\Z$--gradation
on cohomology coincides with the one induced by the principal
gradation on $\ab_+$. Therefore the space $I(\G)$ is linearly spanned
by elements $I_m, m \in I$, where $\deg I_m=m$.

The proof of \lemref{ide} is given in the Appendix. Here we only
explain how to construct an integral of motion in the affine Toda
field theory starting from a class in the first cohomology of the
complex $F^\bu(\G)$.

Consider such a class $\wt{H} \in \oplus_{0\leq i\leq \el}
\pi_{-\al_i}$. Since $\pa = p_{-1}^R \in \ab_-^R$ acts trivially on
cohomologies of the complex $\F^\bu(\G)$ (see , $\pa \wt{H}$ is a
coboundary, i.e. there exists $H \in \pi_0$ that $\bar{\delta}^0 \cdot
H = \pa \wt{H}$.

By construction, the element $H$ has the property that $\q_i \cdot H
\in \pi_{-\al_i}$ is a total derivative for $i = 0,\ldots,\ell$. But
it itself is not a total derivative, because otherwise $\wt{H}$ would
also be a trivial cocycle. Therefore, $\int H dz \neq 0$. But then
$\int H dz$ is an integral of motion of the mKdV hierarchy, because by
construction $\bar{\delta}^0 \cdot \int H dz = \int (\bar{\delta}^0
\cdot H) \; dz = 0$ and hence $\q_i \cdot \int H dz = 0$ for any
$i=0,\ldots,\ell$.

We denote by $H_m$ the integral of motion corresponding to $p_{-m} \in
\ab_-$ via \lemref{ide}. By \lemref{ide}, $I_m = \int H_m dz, m \in
I$, span $I(\G)$. Explicit formulas for $\wt{H}_m$ and $H_m$ can be
found in \cite{EF}, Sect.~5.

\subsection{Proof of \thmref{walg}}

Let $\ovl{W}$ be the Weyl group of $\g$. It is a finite subgroup of
the affine Weyl group $W$ of $\G$. Consider the subcomplex $F^\bu(\g)$
of $F^\bu(\G)$, where
$$
F^j(\g) = \oplus_{w \in \ovl{W}} \pi_{w(\rho)-\rho}.
$$
In particular,
$$F^0(\g) = \pi_0, \quad \quad F^1(\G) = \oplus_{i=1}^\el
\pi_{-\al_i},$$ and the differential $F^0(\g) \arr F^1(\g)$ equals
${\mc H}$ given by formula \eqref{mch} with $J=\{ 1,\ldots,\el \}$.

Consider the quotient complex $F^\bu(\g)/(\C \oplus \on{Im} \pa)$. We
have, by definition, $$\F^0(\g) = \F_0, \quad \quad \F^1(\g) =
\oplus_{1\leq i\leq \el} \F_{-\al_i},$$ and the differential $\F^0(\g)
\arr \F^1(\g)$ is given by $\sum_{0\leq i\leq \el} \q_i$. Therefore
the $0$th cohomology of the complex $\F^\bu(\g)$ is isomorphic to the
$I(\g)$.

\begin{lem}[\cite{FF:laws}]
The cohomology of the complex $F^\bu(\g)$ is isomorphic to
$H^\bu(\nn_+,\pi_0)$. Hence the $0$th cohomology of $F^\bu(\g)$ is
isomorphic to $\C[v_i^{(n)}] \subset \C[u_i^{(n)}]$, and the higher
cohomologies vanish.
\end{lem}

\thmref{walg} immediately follows from this lemma. Thus, the space of
integrals of motion of the Toda theory associated to $\g$ is the
Poisson algebra $\W(\g)$, which is a Poisson subalgebra of $\F_0$. It
is called the {\em classical $\W$--algebra} associated to $\g$.

\subsection{KdV hierarchy is hamiltonian}    \label{mkdvham}

Denote by $\eta_m$ the derivation $\xi_{H_m}$, where $H_m$ is the
density of the $m$th integral of motion of the affine Toda field
theory associated to $\G$ (it is defined up to adding a total
derivative). In particular, we can choose as the one-cocycle
$\wt{H}_1$, the vector $\sum_{0\leq i\leq \el} e^{-\phi_i}$. Then $\pa
\wt{H}_1 = - \sum_{0\leq i\leq \el} u_i e^{-\phi_i}$ and $H_1 =
\frac{1}{2} \sum_{1\leq i\leq \el} u_i u^i$, where $u^i$'s are dual to
$u_i$'s with respect to the inner product on $\h$. Hence $\eta_1 =
\pa$, by formula \eqref{pa}.

Now we have for each $m \in I$, a hamiltonian vector field $\eta_m$ on
$\UU$. On the other hand, we have a vector field $p_{-m}^R$ coming
from the right infinitesimal action of the Lie algebra $\ab_- \subset
\G$ on $N_+/A_+ \simeq \UU$; these vector fields define the flows of
the mKdV hierarchy.

\begin{thm}[\cite{FF:kdv}, Theorem 3]    \label{osnovnoi}
{\em The vector field $\eta_m$ coincides with the vector field $p_{-m}^R$
up to a non-zero constant multiple for any $m \in I$.}
\end{thm}

This, we can rescale $H_m$ so as to make $\eta_m=p_{-m}^R$. The
theorem means that the vector field $p_{-m}^R$ is {\em hamiltonian},
with the hamiltonian being the $m$th integral of motion $H_m$ of the
affine Toda field theory. Note that we already know it when $m=1$,
since $\eta_1 = \pa = p_{-1}^R$.

This means that the mKdV equations are hamiltonian, that is the $m$th
equation of the mKdV hierarchy can be written as
$$
\pa_n u_i(z) = \{ u_i(z),\int H_m dz \}
$$
(up to constant multiple). In this formula $u_i(z)$ stands for the
delta-like functional whose value at $u_i: \R \arr \h$ is $u_i(z)$
(see \cite{FF:laws}). Now recall that $\W(\g)$ is a Poisson subalgebra
of $\F_0$, by \thmref{walg}. It is clear from the definition that
$H_m$ can be chosen in such a way that it lies in $\C[v_i^{(n)}]
\subset \C[u_i^{(n)}]$, so that $\int H_m dz \in \W(\g)$. Then we see
that the KdV hierarchy is also hamiltonian, i.e., it can be written as
$$
\pa_n v_i(z) = \{ v_i(z),\int H_m dz \}.
$$

\begin{prop} The integrals of motion of the affine Toda field theory
(equivalently, the hamiltonians of the corresponding (m)KdV hierarchy)
commute with each other: $$\{ \int H_n dz,\int H_m dz \} = 0$$ in
$\F_0$ for any $n, m \in I$.
\end{prop}

\begin{proof} Since $p_{-m}, m \in I,$ lie in a commutative Lie
algebra, they commute with each other; so do the corresponding vector
fields. By \thmref{osnovnoi}, the same holds for the vector fields
$\eta_m, m \in I$: $[\eta_n,\eta_m]=0$. By formula \eqref{com},
injectivity of the map $\xi$ on $\F_0$, and the definiton of the
vector fields $\eta_m$, the corresponding integrals of motion also
commute with each other.
\end{proof}

We conclude that the KdV and mKdV hierarchies are completely
integrable hamiltonian systems.

\subsection{Quantization}

In conclusion of this lecture, I would like to describe the problem of
quantization of Toda integrals of motion, which was the original
motivation for developing the formalism of these lectures.

The map $\xi: \pi_0 \arr \on{End} \pi_0$ defined in Sect.~5.2 can be
quantized in the following sense. Let $\pi_\la^\hbar =
\pi_\la[[\hbar]]$. There exists a linear map $\xi^\hbar: \pi_0^\hbar
\arr \on{End} \pi_0^\hbar$ such that

\begin{enumerate}

\item[(1)] $\xi^\hbar$ factors through $\F_0^\hbar = \F_0[[\hbar]] =
\pi_0^\hbar/(\on{Im} \pa \oplus \C)[[\hbar]]$;

\item[(2)] the bracket $[ \cdot,\cdot ]_\hbar: \F_0^\hbar \times
\F_0^\hbar \arr \F_0^\hbar$, defined by the formula $$[ \int P dz,\int
R dz ]_\hbar = \int (\xi_P^\hbar \cdot R) \; dz,$$ where $\int$
denotes the projection $\pi_0^\hbar \arr \F_0^\hbar$, is a Lie
bracket;

\item[(3)] $\xi^\hbar = \hbar \xi^{(1)} + \hbar^2(\ldots)$, and the
map $\pi_0 \arr \on{End} \pi_0$ induced by $\xi^{(1)}$ coincides with
$\xi$.

\end{enumerate}

Such a map $\xi^\hbar$ can be defined using the vertex operator
algebra structure on $\pi_0$. This is explained in detail in
\cite{FF:laws}, Sect.~4 (where $\hbar$ is denoted by $\beta^2$).  The
key observation that enables us to construct this map is that $\pi_0$
can be viewed as a Fock representation of a Heisenberg algebra with
generators $b_i(n), i=1,\ldots,\el; n \in \Z$, and relations
$$[b_i(n),b_j(m)] = \hbar n (\al_i,\al_j) \delta_{n,-m}.$$ This module
is endowed with a canonical vertex operator algebra structure
$Y^\hbar: \pi^\hbar_0 \arr \on{End} \pi^\hbar_0[[z,z^{-1}]], P \arr
Y^\hbar(P,z)$. Given $P \in \pi_0^\hbar$, we denote by $\xi_P^\hbar$
the linear endomorphism of $\pi_0^\hbar$ given by the residue,
i.e. the $(-1)$st Fourier component, of $Y^\hbar(P,z)$. This gives us
a map $\xi^\hbar: \pi_0^\hbar \arr \on{End} \pi_0^\hbar$, which
satisfies the conditions above.

Thus, the Gelfand-Dickey-Dorfman structure on $\pi_0$ can be viewed as
a classical limit of the structure of vertex operator algebra on
$\pi^\hbar_0$.

In \cite{FF:laws} we also defined quantum deformations of the maps
$\xi: \pi_0 \arr \on{End} \pi_\la$ and $\pi_\la \arr
\on{Hom}(\pi_0,\pi_\la)$. This enabled us to quantize the operators
$\q_i: \pi_0 \arr \pi_{-\al_i}$ and $\Q_i: \F_0 \arr
\F_{-\al_i}$. Hence we can define the space of quantum integrals of
motion as the intersection of the kernels of the quantum operators
$\Q_0^\hbar,\ldots,\Q_\el^\hbar$. This space could {\em a priori} be
``smaller'' than the space $I(\G)$ of classical integrals of motion,
i.e., it could be that some (or even all) of them do not survive
quantization. However, we proved in \cite{FF:laws} that all integrals
of motion of affine Toda field theory can be quantized.

Our proof was based on the fact that the quantum operators
$\q_i^\hbar$ in a certain sense generate the quantized universal
enveloping algebra $U_q(\n_+)$, where $q=\exp(\pi i \hbar)$ (recall
that the operators $\q_i$ generate $U(\n_+)$). Using this fact, we
were able to deform the whole complex $F^\bu(\G)$ and derive the
quantization property from a deformation theory argument (see
\cite{FF:laws}).

\section*{Appendix}

\subsection*{A.1} {\bf Proof of \propref{bk}.} The proposition has
been proved by B.~Kostant (private communication). The proof given
below is different, but it uses the ideas of Kostant's proof.

Recall that for the non-twisted $\G=\g \otimes \C\ti$,
$$p_{-1} = \pp + f_0,$$ where $\ovl{p}_{-1}$ is given by formula
\eqref{ovlp} and $f_0 = e_{\on{max}} \otimes t^{-1}$. Here
$e_{\on{max}}$ is a generator of the one-dimensional subspace of
$\ovl{\n}_+$ corresponding to the maximal root. More generally, for
$i=1,\ldots,\el$, we can also write
\begin{equation}    \label{pi}
p_{-d_i} = \ovl{p}_{-i} \otimes 1 + r_i \otimes t^{-1},
\end{equation}
where $\ovl{p}_{-i} \in \ovl{\n}_-$, and $r_i \in \ovl{\n}_+$. It is
clear that $[p_{-d_i},p_{-d_j}]=0$ implies
\begin{equation}    \label{allcomm}
[\ovl{p}_{-i},\ovl{p}_{-j}] = [r_i,r_j] = 0.
\end{equation}

Since $[x,e_{\on{max}}]=0$ for all $x \in \ovl{\n}_+$,
$$[f_0,p_{-d_i}] = -[\ovl{p}_{-1},p_{-d_i}] = -
[\ovl{p}_{-1},\ovl{p}_{-i} + r_i \otimes t^{-1}]$$ $$= -
[\ovl{p}_{-1},r_i] \otimes t^{-1}.$$ Furthermore, we find from formula
\eqref{gai}:
$$\ga_i = - (\on{ad} \ovl{p}_{-1})^{h-d_i} \cdot r_i \otimes t^{-1}
\in \h \otimes t^{-1}.$$ Hence the linear independence of the vectors
$\ga_i$ is equivalent to the linear independence of the vectors
$$\ovl{\ga}_i = (\on{ad} \ovl{p}_{-1})^{h-d_i} \cdot r_i, \quad \quad
i=1,\ldots,\el.$$

Let $\{ e,h,f \}$ be the principal $\sw_2$ subalgebra of $\g$,
such that $f=\ovl{p}_{-1}$ and $h=2\ovl{\rho}^\vee$. Recall from
\cite{Ko} that as a principal $\sw_2$--module, $\g$ splits into
the direct sum of irreducible representations $R_i$ of dimension
$2d_i+1$, where $i=1,\ldots,\el$. The multiplicity of $R_i$ in the
decomposition of $\g$ equals the multiplicity of the exponent
$d_i$. Note that different components $R_j$ are mutually orthogonal
with respect to the invariant inner product $(\cdot,\cdot)$ on $\g$,
and hence every element of $\g$ can be written canonically as a sum of
its projections on various $R_j$'s.

The linear independence of the vectors $\ovl{\ga}_i$ is equivalent to
the statement that each $r_i$ has a non-zero projection on $R_i$ (and
moreover, if $d_i$ has multiplicity two, then the projections of the
corresponding $r_i^1, r_i^2$ on $R_i$ are linearly independent). Note
that $\deg r_i = \el-d_i$, with respect to the gradation operator
$\ovl{\rho}^\vee$. Thus, it suffices to show that
$(r_i,\ovl{p}_{-\el+i}) \neq 0$ (resp., the pairing between
$\on{span}(r_i^1,r_i^2)$ and
$\on{span}(\ovl{p}_{-\el+i}^1,\ovl{p}_{-\el+i}^2)$ induced by
$(\cdot,\cdot)$ is non-degenerate).

Now recall that according to Kac \cite{Kac}, Lemma 14.4, the inverse
image of $\ab$ in $\GG$ is a (non-degenerate) Heisenberg Lie
subalgebra $\ab \oplus \C C$. Thus, $[p_n,p_{-n}] \neq 0, \forall n
\in I$ (and moreover, if $n$ has multiplicity two, the pairing between
$\on{span}(p_n^1,p_n^2)$ and $\on{span}(p_{-n}^1,p_{-n}^2)$ induced by
the commutator is non-degenerate). But note that
\begin{equation}    \label{pi1}
p_{d_i} = \ovl{p}_{-\el+i} \otimes t + r_{\el-i} \otimes 1.
\end{equation}
We find from formulas \eqref{pi} and \eqref{pi1} the following
commutator in $\GG$:
$$
[p_{d_i}^k,p_{-d_i}^l] = (\ovl{p}_{-\el+i}^k,r_i^l) C.
$$
The proposition follows.

\subsection*{A.2. The BGG resolution.}
Vector $y$ from the Verma module $M_\la$ is called {\em singular vector} of
weight $\mu \in \h^*$, if $\n_+ \cdot y = 0$ and $x \cdot y = \mu(x) y$ for
any $x \in \h$. We have $M_\la \simeq U(\n_-) \cdot v_\la$, where $v_\la$,
which is called the highest weight vector, is a generator of the space
$\C_\la$.  This vector is a singular vector of weight $\la$. Any singular
vector of $M_\la$ of weight $\mu$ can be uniquely represented as $P \cdot
v_\la$ for some element $P \in U(\n_-)$ of weight $\mu - \la$. This
singular vector canonically defines a homomorphism of $\G$-modules $i_P:
M_\mu \arr M_\la$, which sends $u \cdot v_\mu$ to $(uP) \cdot v_\la$ for
any $u \in U(\n_-)$.  Denote by $i_P^*$ the dual homomorphism $M^*_\la \arr
M^*_\mu$.

There is an isomorphism $U(\n_-) \arr U(\n_+)$, which maps the
generators $f_0,\ldots,f_\el$ to $-e_0,\ldots,-e_\el$. Denote by
$\ovl{P}$ the image of $P \in U(\n_-)$ under this isomorphism.

The homomorphism sending $\al \in \n_+$ to $\al^L$, can be extended in
a unique way to a homomorphism from $U(\n_+)$ to the algebra of
differential operators on $N_+$. Denote the image of $u \in U(\n_+)$
under this homomorphism by $u^L$.

\bigskip

\noindent {\bf Proposition A.1} (\cite{FF:kdv}, Prop.~2){\bf .}  {\em
If $P \cdot v_\la$ is a singular vector in $M_\la$ of weight $\mu$,
then the homomorphism $i^*_P: M^*_\la \arr M^*_\mu$ is given by the
differential operator $\ovl{P}^L$.}

\bigskip

In the case of simple Lie algebras, these homomorphisms were first
studied by B.~Kostant \cite{Kos}. Using Proposition A.1, one can
explicitly construct the differentials of the dual BGG resolution
$B^\bu(\G)$ of \secref{bggres}.

It is known that for each pair of elements of the Weyl group, such
that $w \prec w'$, there is a singular vector $P_{w,w'} \cdot
v_{w(\rho) - \rho}$ in $M_{w(\rho) - \rho}$ of weight $w'(\rho) -
\rho$. By Proposition A.1, this vector defines the homomorphism
$\ovl{P}_{w,w'}^L: M^*_{w(\rho) - \rho} \arr M^*_{w'(\rho) - \rho}$.

It is possible to normalize all $P_{w,w'}$'s in such a way that
$P_{w'_1,w''} P_{w,w_1'} = P_{w_2',w''} P_{w,w_2'}$ for any quadruple of
elements of the Weyl group, satisfying $w \prec w_1', w_2'
\prec w''$. Then we obtain: $\ovl{P}^L_{w,w_1'} \ovl{P}^L_{w'_1,w''} =
\ovl{P}^L_{w,w_2'} \ovl{P}^L_{w_2',w''}$.

The differential $\delta^j: B^j(\G) \arr B^{j+1}(\G)$ of the BGG
complex can be written as follows: $$\delta^j =
\sum_{l(w)=j,l(w')=j+1,w\prec w'} \epsilon_{w,w'} \;
\ovl{P}_{w,w'}^L,$$ where $\epsilon_{w,w'} = \pm 1$ are chosen as in
\cite{bgg,rocha}.

\subsection*{A.3. Proof of \lemref{ide}.}
Since $\pa = p_{-1}^R$ commutes with the differentials of the complex
$F^\bu(\G)$, we can consider the double complex
\begin{equation}    \label{double} \tag{A.3}
\C \larr F^\bu(\G) \stackrel{\pm p_{-1}}{\larr} F^\bu(\G) \larr \C.
\end{equation}
Here $\C \arr \pi_0 \subset F^\bu(\G)$ and $F^\bu(\G) \arr \pi_0 \arr \C$ are
the embedding of constants and the projection on constants, respectively.
We place $\C$ in dimensions $-1$ and $2$ of our complex, and $F^\bu(\G)$ in
dimensions $0$ and $1$.

In the spectral sequence, in which $\pm p_{-1}$ is the $0$th differential,
the first term is the complex $\F^\bu(\G)[-1]$, where $$\F^j(\G) \simeq
\oplus_{l(w)=j} \F_{w(\rho) - \rho}.$$ Indeed, if $\la \neq 0$, then in the
complex $$\pi_\la \stackrel{p_{-1}}{\larr} \pi_\la$$ the $0$th cohomology is
$0$, and the first cohomology is, by definition, the space $\F_\la$. If
$\la = 0$, then in the complex $$\C \larr \pi_0 \stackrel{p_{-1}}{\larr}
\pi_0 \larr \C$$ the $0$th cohomology is $0$ and the first cohomology is,
by definition, the space $\F_0$.

We need to find the $0$th cohomology of the complex $\F^\bu(\G)$, which
is the same as the $1$st cohomology of the double complex
\eqref{double}.

We can compute this cohomology, using the other spectral sequence
associated to our double complex. Since $H^\bu(F^\bu(\G)) \simeq
\bigwedge^\bu(\ab_-)$, we obtain in the first term the following complex
$$\C \larr \wedge^\bu(\ab_-) \stackrel{\pm p_{-1}}{\larr}
\wedge^\bu(\ab_-) \larr \C.$$ By \lemref{action}, the action of
$p_{-1}$ on $\bigwedge^\bu(\ab_-)$ is trivial and hence the
cohomology of the double complex \eqref{double} is isomorphic to
$\bigwedge^\bu(\ab_-)/\C \oplus \bigwedge^\bu(\ab_-)/\C[-1]$. In
particular, we see that the space of Toda integrals is isomorphic to
$\ab_-$.

\end{document}